\documentclass[11pt]{article}
\usepackage[utf8]{inputenc}
\usepackage[T1]{fontenc}
\usepackage[frenchb,english]{babel} 
\usepackage{textcomp}
\usepackage{amsmath,amssymb}
\usepackage{amsthm}
                  
\usepackage{lmodern}

\usepackage[a4paper,margin=2cm]{geometry}

\usepackage{graphicx}             
\usepackage{xcolor}               
\usepackage{microtype,stmaryrd}          

\usepackage{hyperref}
\hypersetup{pdfstartview=XYZ}

\usepackage{bbm}
\usepackage{mathrsfs}
\usepackage{subfig}

\def\build#1_#2^#3{\mathrel{
\mathop{\kern 0pt#1}\limits_{#2}^{#3}}}
\def\llbracket{[\hspace{-.10em} [ }
\def\rrbracket{ ] \hspace{-.10em}]}

\newtheorem{theorem}{Theorem}
\newtheorem{proposition}[theorem]{Proposition}
\newtheorem{definition}[theorem]{Definition}
\newtheorem{lemma}[theorem]{Lemma}

\def\w{\mathrm{w}}
\def\t{\mathcal{T}}
\def\f{\mathcal{F}}
\def\r{\mathcal{R}}

\def\W{\mathcal{W}}
\def\S{\mathcal{S}}

\def\N{\mathbb{N}}
\def\M{\mathbb{M}}

\def\D{\mathbb{D}}
\def\P{\mathbb{P}}
\def\E{\mathbb{E}}

\def\R{\mathbb{R}}

\def\ee{\mathcal{E}}
\def\ve{{\varepsilon}}
\def\la{\longrightarrow}

\def\dd{\mathrm{d}}
\def\wh{\widehat}
\def\wt{\widetilde}

\def\nn{\mathcal{N}}

\def\H{\mathbb{H}}

\def\be{\mathbf{e}}

\def\PP{\mathbf{P}}

\def\rem{\noindent{\bf Remark. }}

\title{The Brownian disk viewed from a boundary point\footnote{Supported by the ERC Advanced Grant 740943 {\sc GeoBrown}}}

\author{Jean-Fran\c cois Le Gall}

\date{\small Universit\'e Paris-Saclay}
\begin{document}
\maketitle

\begin{abstract}
We provide a new construction of Brownian disks in terms of forests of
continuous random trees equipped with nonnegative labels corresponding to distances 
from a distinguished point uniformly distributed on the boundary of the disk. 
This construction shows in particular that distances from the distinguished point
evolve along the boundary as a five-dimensional Bessel bridge. As 
an important ingredient of our proofs, we show that the uniform measure
on the boundary, as defined in the earlier work of Bettinelli and Miermont, is the limit 
of the suitably normalized volume measure on a small tubular neighborhood 
of the boundary. 
Our construction also yields
a simple proof of the equivalence between the two definitions of the 
Brownian half-plane.
\end{abstract}

\tableofcontents

\section{Introduction}

Brownian disks are random compact metric spaces that serve as models of random
geometry and arise as scaling 
limits of large random planar maps with a boundary \cite{AHS,BMR,Bet,BM,GM2}. They appear 
as special subsets of the Brownian map, and in particular as connected components of the 
complement of balls in the Brownian map \cite{Disks}. Brownian disks are also closely related
to the Liouville quantum gravity surfaces called quantum disks, see \cite[Corollary 1.5]{MS}, as well as the survey \cite{Miller} and the
references therein. The initial construction of Brownian disks was given by Bettinelli \cite{Bet} 
in terms of a forest of continuous random trees
equipped with Brownian labels. In this construction, labels correspond to distances from a distinguished 
point belonging to the interior of the Brownian disk. A different construction still based on a labeled continuous
random tree appeared in \cite{Disks}, with labels now corresponding to distances from the boundary of
the disk. The main goal of the present work is to present a new construction of Brownian disks where labels represent
distances from a point chosen uniformly at random on the boundary. This construction has several
interesting consequences. In particular, it shows that, if one starts from a distinguished point chosen at random on the
boundary and then moves along the boundary, distances 
from the distinguished point evolve exactly like a five-dimensional 
Bessel bridge. In contrast with preceding constructions \cite{Bet,BM,Disks}, we do
not rely on discrete approximations to establish the validity of our method, but
rather we pass to the limit in the construction of \cite{Bet} by letting the distinguished 
point tend to the boundary.

Let us give an informal description of our construction (see Section \ref{sec:dis-bdry} for a more precise 
presentation). We start
from the circle, which we view as the interval $[0,1]$ with the two points $0$ and $1$ identified, 
and we assign  a ``label'' $\Lambda_t$ to each point $t$ of $[0,1]$, in such a way
that the process $(\Lambda_t)_{t\in[0,1]}$ is a five-dimensional 
Bessel bridge from $0$ to $0$ scaled by the factor $\sqrt{3}$. 
We then consider a Poisson forest of continuous 
random trees (scaled versions of the celebrated Aldous Brownian CRT) that are rooted randomly
on the circle. For every tree $\t$ in this forest, and for every point $u\in \t$, we assign a label $\Lambda_u$
to $u$: The collection $(\Lambda_u)_{u\in\t}$ is distributed as Brownian motion 
indexed by $\t$, started from 
the label of the root (recall that the root of $\t$ belongs to the circle).
We let $\mathfrak{H}$ be the geodesic metric space consisting of the union of the circle and 
the collection of those trees that have only nonnegative labels (we just remove those trees where negative labels occur). In this way, every
point $u$ of $\mathfrak{H}$ has been assigned a nonnegative label $\Lambda_u$.
For $u,v\in\mathfrak{H}$, we set 
$$\Delta^\circ(u,v) = \Lambda_u+\Lambda_v-2\,\max\Big(\min\{\Lambda_w:w\in [|u,v|]\},\min\{\Lambda_w:w\in [|v,u|]\}\Big),$$
where 
$[|u,v|]$ is the ``interval'' of $\mathfrak{H}$ consisting of points visited when going from $u$ to $v$ in 
 ``clockwise direction'' along $\mathfrak{H}$ (see Section \ref{sec:BD-BHP} for more precise definitions). Finally, we define
 $\Delta(u,v)$ as the largest pseudo-metric on $\mathfrak{H}$ that is bounded above by $\Delta^\circ(u,v)$. 
 Then
we consider the quotient space $\D:=\mathfrak{H}/\simeq$ for the equivalence relation
defined by setting $x\simeq y$ if and only $\Delta(x,y)=0$. Theorem \ref{cons-disk}
below states that $\D$ equipped with the distance induced by $\Delta$ is a free Brownian disk of
perimeter $1$ pointed
at a uniform boundary point.

The preceding definitions of $\Delta^\circ$ and $\Delta$ are of course very similar to the construction of the 
Brownian map (see e.g.~\cite{Uniqueness}), of the Brownian disk \cite{Bet}, or of the Brownian plane \cite{Plane}.
Indeed, we derive Theorem \ref{cons-disk} by a suitable passage to the limit from the construction of 
the free pointed Brownian disk that is given in
\cite{Bet}  --- note that \cite{Bet} considers the slightly different model of the Brownian disk
with prescribed volume and perimeter, but the same method applies to the free Brownian disk with
minor changes. In this construction, labels correspond to distances from a distinguished point
distributed according to the volume measure of the Brownian disk (see formula \eqref{pointed-no}
below for a more precise statement describing the distribution of the distinguished point). The idea is then to condition the distinguished 
point to lie within distance at most $\ve$ from the boundary $\partial\D$ and to pass to the limit $\ve\to0$. For this 
passage to the limit, it is crucial to have information about the probability measure $\mu_\ve$ obtained 
by normalizing the restriction of the volume measure of $\D$ to the tubular neighborhood of radius $\ve$ of the boundary.
More precisely, one needs the fact that $\mu_\ve$
converges when $\ve \to 0$ to the uniform probability measure $\mu$ on the boundary, as defined in the construction
of \cite{Bet,BM}. The convergence of $\mu_\ve$ as $\ve\to 0$ towards a 
probability measure $\nu$ supported on $\partial \D$ had already been
obtained in \cite{Disks}, but the equality $\mu=\nu$ was still open. 
Theorem \ref{approx-unif-bdry} below shows that this equality holds, so that 
the two natural ways of defining a uniform measure on the boundary are indeed equivalent.

Another important ingredient consists in studying the behavior of labels on the boundary, under the condition 
that the distinguished 
point lies within distance at most $\ve$ from $\partial\D$. In the construction of \cite{Bet,BM}, labels
evolve along the boundary like a Brownian bridge scaled by the constant $\sqrt{3}$, and one may replace the Brownian bridge by
a normalized Brownian excursion $\be=(\be_t)_{0\leq t\leq 1}$ thanks to Vervaat's transformation \cite{Vervaat}.
Under the preceding conditioning, labels along the boundary evolve like $\sqrt{3}\,\be^\ve$, where 
the distribution of $\be^\ve$ is specified by
\begin{equation}
\label{condi-exc}
\E[F(\be^\ve)]:= C_\ve^{-1}\,\E\Big[ F(\be)\,\exp\Big(-\int_0^1 \frac{\dd t}{(\ve + \be_t)^2}\Big)\Big],
\end{equation}
where $C_\ve$ is the appropriate normalizing constant. Proposition \ref{conv-bridge} below
states that $\be^\ve$ converges in distribution as $\ve\to 0$ to a five-dimensional Bessel bridge $\mathbf{b}=(\mathbf{b}_t)_{0\leq t\leq 1}$.
For our applications, we need in fact a more precise result showing that, for every $\delta\in(0,1/2)$, it is possible 
to couple $\be^\ve$ and $\mathbf{b}$ so that the equality $\ve + \be^\ve_t=\mathbf{b}_t$ holds
for every $t\in[\delta,1-\delta]$, with high probability when $\ve\to 0$.

Our construction of the Brownian bridge is closely related to the definition of the Brownian half-plane 
proposed by Caraceni and Curien \cite{CC}, which involves a two-sided five-dimensional Bessel process.
Another definition of the Brownian half-plane, which is close to the Bettinelli construction of Brownian
disks, has been given independently by Gwynne and Miller \cite{GM0} and by Baur, Miermont and Ray \cite{BMR},
but initially it was not clear that this definition yields the same random object as the Caraceni-Curien definition
(see the comments in \cite[Remark 2.7]{BMR} and in \cite[Section 1.6]{GM0}). 
Recently, Budzinski and Riera \cite{BR} have been able to prove the equivalence of the two definitions via 
discrete approximations. In Section \ref{sec:CC}, we provide a short simple proof of this equivalence based
on our new construction of the Brownian disk. 

The paper is organized as follows. Section \ref{prelim} contains a few preliminaries. In particular, we 
recall the formalism of snake trajectories, which provides a convenient framework to deal
with continuous random trees equipped with labels, and we define 
the spaces of compact or non-compact measure metric spaces that are relevant to the present work.
The technical Section \ref{conv-Bessel} investigates  the limiting behavior of the ``excursions'' $\be^\ve$
distributed as in \eqref{condi-exc}. We start by recalling several properties of Bessel processes, and 
especially of first-passage Bessel bridges, as these properties play an important role in the proof of the
key technical Proposition \ref{conv-bridge}. Section \ref{sec:BD-BHP} is mainly devoted to
recalling the constructions
of the (free pointed) Brownian disk and of the Brownian half-plane found in \cite{BMR,Bet,BM,GM0}.
Our presentation is slightly different from the latter papers and adapted to our purposes. 
In Section \ref{sec:approx}, we prove Theorem \ref{approx-unif-bdry}  concerning the approximation
of the uniform measure on the boundary by the volume measure on a tubular neighborhood of small radius.
Here the Brownian half-plane is used as a tool: We derive a Brownian half-plane version of the desired approximation
via an application of the ergodic theorem, and we then
use a suitable coupling of the Brownian disk and the Brownian half-plane near a 
boundary point. In Section \ref{sec:dis-bdry}, we prove Theorem \ref{cons-disk} giving our
 construction of the free Brownian disk pointed at a uniform boundary point. Here, the method 
 consists in coupling the Brownian disk pointed at a point lying within distance $\ve$ from the 
 boundary, and the candidate space for the Brownian disk pointed at a uniform boundary point, in such a way that
 one can get a suitable bound on
 the Gromov-Hausdorff distance between these two spaces. For this coupling, the precise
 statement of Proposition \ref{conv-bridge} is crucial. Finally, Section \ref{sec:CC} proves the
 equivalence of the two definitions of the Brownian half-plane.
 
 \medskip
 
 \noindent{\bf Acknowledgements.} I thank Nicolas Curien and Armand Riera for stimulating conversations. 


\section{Preliminaries}
\label{prelim}

\subsection{Snake trajectories}
\label{sec:preli}

We use the formalism of snake trajectories to deal with continuous random trees whose
vertices are assigned real labels. In this section, we briefly recall the 
notation and definitions that are relevant to the present work.
We refer to \cite{ALG} for more details. 

We denote the space of all finite (real) paths by $\W$. Here a finite 
path $\w$ is a continuous mapping $\w:[0,\zeta]\la\R$, where the
number $\zeta=\zeta_{(\w)}\geq 0$ is called the lifetime of $\w$. 
The endpoint or tip of the path $\w$ is denoted by $\wh \w=\w(\zeta_{(\w)})$.
The space
$\W$ is equipped with the distance
$$d_\W(\w,\w')=|\zeta_{(\w)}-\zeta_{(\w')}|+\sup_{t\geq 0}|\w(t\wedge
\zeta_{(\w)})-\w'(t\wedge\zeta_{(\w')})|.$$
We set $\W_0:=\{\w\in\W:\w(0)=0\}$. The trivial path of $\W_0$ with zero lifetime is identified to the point $0$ of $\R$.

\begin{definition}
\label{def:snakepaths}
A snake trajectory $\omega$ (with initial point $0$) is a continuous mapping $s\mapsto \omega_s$
from $\R_+$ into $\W_0$ 
which satisfies the following two properties:
\begin{enumerate}
\item[\rm(i)] We have $\omega_0=0$ and the number $\sigma(\omega):=\sup\{s\geq 0: \omega_s\not =0\}$,
called the duration of the snake trajectory $\omega$,
is finite (by convention $\sup\varnothing=0$). 
\item[\rm(ii)] {\rm (Snake property)} For every $0\leq s\leq s'$, we have
$\omega_s(t)=\omega_{s'}(t)$ for every $t\in[0,\displaystyle{\min_{s\leq r\leq s'}} \zeta_{(\omega_r)}]$.
\end{enumerate} 
\end{definition}

We denote the set of all snake trajectories by $\S$. 
If $\omega\in \S$, we often write $W_s(\omega)=\omega_s$ and $\zeta_s(\omega)=\zeta_{(\omega_s)}$
for every $s\geq 0$. The set $\S$ is equipped with the distance
$$d_{\S}(\omega,\omega')= |\sigma(\omega)-\sigma(\omega')|+ \sup_{s\geq 0} \,d_\W(W_s(\omega),W_{s}(\omega')).$$
It is not hard to verify that a snake trajectory $\omega$ is completely determined 
by the knowledge of the lifetime function $s\mapsto \zeta_s(\omega)$ and of the tip function $s\mapsto \wh W_s(\omega)$ (see \cite[Proposition 8]{ALG}).

Let $\omega\in \S$ be a snake trajectory. The lifetime function $s\mapsto \zeta_s(\omega)$ codes a
compact $\R$-tree, which will be denoted 
by $\t_{(\omega)}$. This $\R$-tree is the quotient space $\t_{(\omega)}:=[0,\sigma(\omega)]/\!\sim$ 
of the interval $[0,\sigma(\omega)]$
for the equivalence relation
$$s\sim s'\ \hbox{if and only if }\ \zeta_s(\omega)=\zeta_{s'}(\omega)= \min_{s\wedge s'\leq r\leq s\vee s'} \zeta_r(\omega),$$
and $\t_{(\omega)}$ is equipped with the distance induced by
$$d_{(\omega)}(s,s')= \zeta_s(\omega)+\zeta_{s'}(\omega)-2 \min_{s\wedge s'\leq r\leq s\vee s'} \zeta_r(\omega).$$
(see e.g.~\cite[Section 3]{LGM} for more information about the
coding of $\R$-trees by continuous functions).  We write $p_{(\omega)}:[0,\sigma(\omega)]\la \t_{(\omega)}$
for the canonical projection. By convention, $\t_{(\omega)}$ is rooted at the point
$\rho_{(\omega)}:=p_{(\omega)}(0)$, and the volume measure on $\t_{(\omega)}$ is defined as the pushforward of
Lebesgue measure on $[0,\sigma(\omega)]$ under $p_{(\omega)}$. 

By property (ii) in the definition of  a snake trajectory, the condition $p_{(\omega)}(s)=p_{(\omega)}(s')$ implies that 
$W_s(\omega)=W_{s'}(\omega)$. So the mapping $s\mapsto W_s(\omega)$ can be viewed as defined on the quotient space $\t_{(\omega)}$.
For $u\in\t_{(\omega)}$, we set $\ell_u(\omega):=\wh W_s(\omega)$, for any $s\in[0,\sigma(\omega)]$ such that $u=p_{(\omega)}(s)$. We interpret $\ell_u(\omega)$ as a ``label'' assigned to the ``vertex'' $u$ of $\t_{(\omega)}$. 
Notice that the mapping $u\mapsto \ell_u(\omega)$ is continuous on $\t_{(\omega)}$. 
We set $W_*(\omega):=\min\{\ell_u(\omega):u\in\t_{(\omega)}\}$ and $W^*(\omega):=\max\{\ell_u(\omega):u\in\t_{(\omega)}\}$,
for $\omega\in\S$.

We now introduce a $\sigma$-finite measure on $\S$ that plays an important role
in the present work.

\begin{definition}
\label{DefN0}
The Brownian snake excursion 
measure $\N_0$ is the $\sigma$-finite measure on $\S$ that is characterized by the following two properties: 
\begin{enumerate}
\item[\rm(i)] The distribution of the lifetime function $(\zeta_s)_{s\geq 0}$ under $\N_0$ is the It\^o 
measure of positive excursions of linear Brownian motion, normalized so that, for every $\ve>0$,
$$\N_0\Big(\sup_{s\geq 0} \zeta_s >\ve\Big)=\frac{1}{2\ve}.$$
\item[\rm(ii)] Under $\N_0$ and conditionally on $(\zeta_s)_{s\geq 0}$, the tip function $(\wh W_s)_{s\geq 0}$ is
a centered Gaussian process with covariance function 
$$K(s,s'):= \min_{s\wedge s'\leq r\leq s\vee s'} \zeta_r.$$
\end{enumerate}
\end{definition}

Informally, property (ii) says that, under $\N_0$ and conditionally on $(\zeta_s)_{s\geq 0}$, the labels $(\ell_u)_{u\in\t_{(\omega)}}$ are distributed as Brownian motion indexed by $\t_{(\omega)}$. The measure $\N_0$ can be interpreted as the excursion measure away from $0$ for the 
Markov process in $\W$ called the Brownian snake
(we refer to 
\cite{Zurich} for a detailed study of the Brownian snake and its excursion measures). 
For our purposes, it will be important to know the distribution of 
the minimum $W_*$ under $\N_0$:
For every $y<0$, we have
\begin{equation}
\label{hittingpro}
\N_0(W_*\leq y)={\displaystyle \frac{3}{2y^2}}.
\end{equation}
See e.g. \cite[Section VI.1]{Zurich} for a proof.

The following scaling property is often useful. For $\lambda>0$, for every 
$\omega\in \S$, we define $\Theta_\lambda(\omega)\in \S$
by $\Theta_\lambda(\omega):=\omega'$, with
$$\omega'_s(t):= \sqrt{\lambda}\,\omega_{s/\lambda^2}(t/\lambda)\;,\quad
\hbox{for } s\geq 0\hbox{ and }0\leq t\leq \zeta'_s:=\lambda\zeta_{s/\lambda^2}.$$
Then it is a simple exercise to verify that $\Theta_\lambda(\N_0)= \lambda\, \N_{0}$. 

Let us introduce some additional notation. 
If $I$ is an interval of $\R$ and $E$ is a metric space, we write $C(I,E)$ for the space
of all continuous functions from $I$ into $E$ (in particular, $C([0,t],\R)$ is
a subset of $\W$, for every $t\geq 0$). 
When $I$ is a compact interval and $E=\R$ or $\R_+$,
the space $C(I,E)$ will be equipped with the topology of uniform convergence, 
and then convergence of probability measures on $C(I,E)$ will be in the
usual sense of weak convergence of probability measures on a Polish space.

We also write $M_p(I\times \S)$ for the set of all
point measures (countable sums of Dirac masses) on $I\times \S$. As usual, $M_p(I\times \S)$
is equipped with the $\sigma$-field generated by the mappings $\gamma\mapsto \gamma(A)$,
when $A$ varies among the Borel subsets of $I\times \S$. 

\subsection{Spaces of compact and locally compact metric spaces}
\label{sec:GHP}

Recall that a compact measure metric space is a compact metric space $(X,d)$ equipped with a Borel
finite measure $\mu$ on $X$, which is sometimes called the volume measure on $X$. 
If there is a distinguished point $x\in X$, we say that $(X,d,\mu,x)$ is a pointed compact measure metric space.

We write $\M$, resp.~$\M^\bullet$, for the set of all compact measure metric spaces, resp.~of all 
pointed compact measure metric spaces, where two such spaces
$(X,d,\mu)$ and $(X',d',\mu')$, resp.~$(X,d,\mu,x)$ and $(X',d',\mu',x')$, are identified if there exists an 
isometry $\phi$ from $X$ onto 
$X'$ such that $\mu'$ is the pushforward of $\mu$ under $\phi$ (and $\phi(x)=x'$ in the pointed case). 
Both $\M$ and $\M^\bullet$ are Polish spaces when equipped with the Gromov-Hausdorff-Prokhorov distance
(see e.g.~\cite[Section 2.1]{Disks} for a definition). 

We will also consider the case of non-compact spaces. We restrict our attention to length spaces
(a metric space $(E,d)$ is called a length space if, for every $x,y\in E$, the distance $d(x,y)$
is the infimum of lengths of continuous paths from $x$ to $y$). Recall also that
a metric space $(E,d)$ is said to be boundedly compact if the closed balls of $E$ are compact.
A length space is boundedly compact if and only if it is locally compact and complete
\cite[Proposition 2.5.22]{BBI}. We let $\M^\bullet_{bcl}$ denote the space of all (isometry classes of)
boundedly compact length spaces $(X,d)$ given with a distinguished point $x$ and 
a measure $\mu$ which is finite on compact subsets of $X$. The set $\M^\bullet_{bcl}$ 
can be equipped with the ``local'' Gromov-Hausdorff-Prokhorov distance as defined
in \cite{ADH} and is then also a Polish space.

\section{Convergence to the five-dimensional Bessel bridge}
\label{conv-Bessel}

The main goal of this section is to prove Proposition \ref{conv-bridge} concerning 
the asymptotic behavior of the ``excursions'' $\be^\ve$ defined in 
\eqref{condi-exc}. Before stating and proving Proposition \ref{conv-bridge}, we need to gather a few facts
about Bessel processes (more information can be found in \cite[Chapter XI]{RY} and especially in \cite{PY2}). It will be convenient to introduce a random process
$R=(R_t)_{t\geq 0}$ and probability measures $\P^{(5)}_x$ and $\P^{(-1)}_x$,
for every $x\geq 0$,
such that $R$ is a five-dimensional Bessel process that starts at $x$ under $\P^{(5)}_x$,
and similarly $R$ is a Bessel process of dimension $-1$ that starts at $x$ under $\P^{(-1)}_x$.
Recall that the Bessel process of dimension $-1$ is absorbed at $0$, and that the five-dimensional Bessel can be viewed as the Bessel process of dimension $-1$ conditioned
to escape to infinity, in the sense of $h$-transforms. More precisely, for every $x>0$ and $t>0$, for every
nonnegative measurable function $F$ on $C([0,t],\R)$, we have
\begin{equation}
\label{h-transform}
\E^{(5)}_x[F((R_s)_{s\leq t})]=\E^{(-1)}_x\Big[ \Big(\frac{R_t}{x}\Big)^{3}F((R_s)_{s\leq t})\Big].
\end{equation}
For every $x\geq 0$, we set $T^{(R)}_x:=\inf\{t\geq 0:R_t=x\}$ and $L^{(R)}_x:=\sup\{t\geq 0:R_t=x\}$
with the usual conventions $\inf\varnothing=\infty$ and $\sup\varnothing=0$.
It follows from \eqref{h-transform} 
that, if $0<\ve<x$, $\P^{(5)}_x(T^{(R)}_\ve<\infty)=(\ve/x)^3$ and
\begin{equation}
\label{h-trans}
(R_t)_{0\leq t\leq T^{(R)}_\ve}\hbox{ under }\P^{(5)}_x(\cdot\mid T^{(R)}_\ve<\infty)\build{=}_{}^{\rm(d)} (R_t)_{0\leq t\leq T^{(R)}_\ve}
\hbox{ under }\P^{(-1)}_x.
\end{equation}
Furthermore, as a consequence of 
Nagasawa's time-reversal theorem \cite[Theorem VII.4.5]{RY}, we have for every $x>0$, 
\begin{equation}
\label{Naga}
(R_t)_{0\leq t\leq L^{(R)}_x}\hbox{ under }\P^{(5)}_0 \build{=}_{}^{\rm(d)} \Big(R_{T^{(R)}_0-t}\Big)_{0\leq t\leq T^{(R)}_0}
\hbox{ under }\P^{(-1)}_x.
\end{equation}
This implies that the process $(L^{(R)}_x)_{x\geq 0}$ has independent increments under $\P^{(5)}_0$. This
property (for more general Bessel processes) was first observed by Getoor \cite{Getoor}.

Fix $0<\ve<x$. By \eqref{Naga}, the law of $T^{(R)}_\ve$ under $\P^{(-1)}_x$ is equal to the
law of $L^{(R)}_x-L^{(R)}_\ve$ under $\P^{(5)}_0$. Thus, for every $\lambda>0$,
$$\E^{(-1)}_x[\exp(-\lambda T^{(R)}_\ve)]= \E^{(5)}_0[\exp(-\lambda(L^{(R)}_x-L^{(R)}_\ve))]=
\frac{\E^{(5)}_0[\exp(-\lambda L^{(R)}_x)]}{\E^{(5)}_0[\exp(-\lambda L^{(R)}_\ve)]}.$$
From the main result of \cite{Getoor}, the density of $L^{(R)}_x$ under $\P^{(5)}_0$ is
the function
$$t\mapsto r_t(x,0):=\frac{x^3}{\sqrt{2\pi t^5}}\,\exp(-\frac{x^2}{2t}),$$
from which one easily computes the Laplace transform 
$\E^{(5)}_0[\exp(-\lambda L^{(R)}_x)]= (1+x\sqrt{2\lambda})e^{-\sqrt{2\lambda}}$. Hence,
$$\E^{(-1)}_x[\exp(-\lambda T^{(R)}_\ve)]= \frac{1+x\sqrt{2\lambda}}{1+\ve\sqrt{2\lambda}}
e^{-(x-\ve)\sqrt{2\lambda}}.$$
We will need the explicit formula for the density of $T^{(R)}_\ve$ under $\P^{(-1)}_x$, which we can obtain 
by inverting the Laplace transform in the preceding display. We note that 
\begin{equation}
\label{Lapl-deco}
\frac{1+x\sqrt{2\lambda}}{1+\ve\sqrt{2\lambda}}e^{-(x-\ve)\sqrt{2\lambda}}= 
e^{-(x-\ve)\sqrt{2\lambda}} + \frac{x-\ve}{\ve} \frac{\ve\sqrt{2\lambda}}{1+\ve\sqrt{2\lambda}} e^{-(x-\ve)\sqrt{2\lambda}}.
\end{equation}
We have $e^{-(x-\ve)\sqrt{2\lambda}} =\int_0^\infty \dd t\,e^{-\lambda t}\,q_t(x,\ve)$, where the function
\begin{equation}
\label{density-hittingBM}
t\mapsto q_t(x,\ve):=\frac{x-\ve}{\sqrt{2\pi t^3}} \exp\Big(-\frac{(x-\ve)^2}{2t}\Big)
\end{equation}
is the density of the hitting time of $\ve$ for a linear Brownian motion started at $x$.
On the other hand, for every $a,b>0$ and $t\geq 0$, set
$$g_{a,b}(t):= e^{ab+a^2t}\,\mathrm{erfc}(a\sqrt{t}+ \frac{b}{2\sqrt{t}}),$$
where we recall the standard notation $\mathrm{erfc}(x)=\frac{2}{\sqrt{\pi}}\int_x^\infty e^{-y^2}\,\dd y$. 
Via an integration by parts and simple calculations, one gets that the Laplace transform of $g_{a,b}$ is
$$\int_0^\infty \dd t\,e^{-\lambda t}\,g_{a,b}(t)= \frac{e^{-b\sqrt{\lambda}}}{\sqrt{\lambda}(a+\sqrt{\lambda})}.$$
Notice that $g_{a,b}(0)=0$ and $g_{a,b}$ tends to $0$ at infinity. It follows that
$$\int_0^\infty \dd t\,e^{-\lambda t}\,g'_{a,b}(t)= \lambda \int_0^\infty \dd t\,e^{-\lambda t}\,g_{a,b}(t) = 
\frac{\sqrt{\lambda}}{a+\sqrt{\lambda}}\, e^{-b\sqrt{\lambda}}.$$
Recalling \eqref{Lapl-deco}, and using the last display with $b=\sqrt{2}(x-\ve)$ and $a=1/(\ve\sqrt{2})$,
we get that the law of $T^{(R)}_\ve$ under $\P^{(-1)}_x$
has a density given by
$$t\mapsto r_t(x,\ve):= q_t(x,\ve) + \frac{x-\ve}{\ve}\,g'_{(\ve\sqrt{2})^{-1},\sqrt{2}(x-\ve)}(t).$$
From the explicit expression for $g_{a,b}$, we get
$$r_t(x,\ve)=\Bigg( \frac{1}{2\ve^3} \mathrm{erfc}\Big(\frac{\sqrt{t}}{\ve\sqrt{2}} + \frac{x-\ve}{\sqrt{2t}}\Big)
\exp\Big( \Big(\frac{\sqrt{t}}{\ve\sqrt{2}} + \frac{x-\ve}{\sqrt{2t}}\Big)^2\Big) 
-\frac{1}{\ve^2\sqrt{2\pi t}}+ \frac{x}{\ve\sqrt{2\pi t^3}}\Bigg)\,(x-\ve)\,e^{-(x-\ve)^2/(2t)}.$$
Using the asymptotic expansion $\mathrm{erfc}(z)\exp(z^2)= (z\sqrt{\pi})^{-1}(1-\frac{1}{2}z^{-2} + o(z^{-2}))$ as $z\to\infty$,
one easily verifies that $r_t(x,\ve)\la r_t(x,0)$ and $r_t(x+\ve,\ve)\la r_t(x,0)$ as $\ve\to 0$. 

We will need to introduce the ``first-passage Bessel bridge'' giving the distribution of
$(R_s)_{0\leq s\leq T^{(R)}_\ve}$ under $\P^{(-1)}_x(\cdot \,|\, T^{(R)}_\ve=t)$, for $0\leq \ve<x$ and $t>0$
(beware that this first-passage bridge should not be confused with the usual Bessel bridges studied in \cite{PY}). Before giving a precise
definition of this bridge, let us  introduce the transition densities of the
Bessel process of dimension $-1$ killed upon hitting $\ve$. For every $\ve\geq 0$, these transition densities are
the (continuous) functions $p^{(\ve)}_t(y,z)$, defined for $t>0$ and $y,z>\ve$ and such that
$$\E^{(-1)}_y[\varphi(R_t)\,\mathbf{1}_{\{T^{(R)}_\ve>t\}}]=\int_{(\ve,\infty)} \dd z\,p^{(\ve)}_t(y,z)\,\varphi(z),$$
for every nonnegative measurable function $\varphi$ on $\R_+$. Let $p_t(y,z)$, $t,y,z>0$, denote the transition densities
of the Bessel process of dimension $-1$. Then, using the strong Markov property at time $T^{(R)}_\ve$, one easily gets,
for $\ve >0$,
\begin{equation}
\label{density-condi}
p^{(\ve)}_t(y,z) = p_t(y,z) - \int_0^t \dd s\,r_s(y,\ve)\,p_{t-s}(\ve,z).
\end{equation}
For $\ve=0$, we have just $p^{(0)}_t(y,z) = p_t(y,z)$. 
Furthermore, we have for every $0<s<t$ and $x>\ve>0$,
\begin{equation}
\label{Markov-densi}
r_t(x,\ve)=\int_\ve^\infty \dd y\,p^{(\ve)}_s(x,y)\,r_{t-s}(y,\ve).
\end{equation}
For $y,z>0$, let $G(y,z)=\int_0^\infty \dd t\,p_t(y,z)$ be the Green function of the Bessel process of dimension $-1$.
Then, $G(y,z)=\frac{2}{3}z(1\wedge \frac{y^3}{z^3})$ (a simple way to get this formula is to
use \eqref{h-transform} to observe that $G(y,z)=\frac{y^3}{z^3}G'(y,z)$, where $G'$ is the Green function of the five-dimensional Bessel process, which
is easily computed from the Green function of Brownian motion). Using \eqref{density-condi}, it follows that, for $y,z>\ve$,
$$G^{(\ve)}(y,z):=\int_0^\infty \dd t\,p^{(\ve)}_t(y,z)= G(y,z)-G(\ve,z)= \frac{2}{3}z(1\wedge \frac{y^3}{z^3})- \frac{2}{3}\frac{\ve^3}
{z^2}.$$
where we made the convention $G(0,z)=0$. If $y>z$, $G^{(\ve)}(y,z)$ does not depend on $y$ and is equal to
$$G^{(\ve)}(\infty,z):=\frac{2}{3}z\,(1-\frac{\ve^3}{z^3}).$$

\begin{proposition}
\label{Bessel-bridge}
Let $x>\ve\geq 0$. For every $t>0$, we can define a probability measure $\Pi^{(x,\ve)}_t(\dd \w)$
on $C([0,t],\R_+)$ in such a way that:
\begin{itemize}
\item[\rm(i)] The collection $(\Pi^{(x,\ve)}_t)_{t> 0}$ is a regular version of the conditional
distributions of $(R_s)_{0\leq s\leq T^{(R)}_\ve}$ knowing $T^{(R)}_\ve=t$ under $\P^{(-1)}_x$.
\item[\rm(ii)] For every $0\leq s<t$, the distribution of $(\w(u))_{0\leq u\leq s}$ under 
$\Pi^{(x,\ve)}_t(\dd \w)$ is absolutely continuous with respect to the distribution 
of $(R_u)_{0\leq u\leq s}$ under $\P^{(-1)}_x$, with a Radon-Nikodym density given by
\begin{equation}
\label{RadonN}
(\w(u))_{0\leq u\leq s} \mapsto \mathbf{1}_{\{\w(u)>\ve:\forall u\in[0,s]\}}\,\frac{r_{t-s}(\w(u),\ve)}{r_t(x,\ve)}.
\end{equation}
\end{itemize}
Furthermore, for every $t>0$, $\Pi^{(x+\ve,\ve)}_t(\dd \w)$ converges weakly to $\Pi^{(x,0)}_t(\dd \w)$ as $\ve\to 0$.
\end{proposition}

In what follows, we will
write $\E^{(-1)}_x[F((R_s)_{0\leq s\leq t})\,|\,T^{(R)}_\ve=t]$ instead of $\int  \Pi^{(x,\ve)}_t(\dd \w)\,F(\w)$
when $F$ is a measurable function on $C([0,t],\R_+)$, . 

\medskip
\rem The proof below applies to the more general setting where the Bessel process of dimension $-1$
is replaced by a Bessel process of dimension $2(1-\nu)$ ($\nu>0$) and the role of the Bessel process 
of dimension $5$ is played by a Bessel process of dimension $2(1+\nu)$. In particular, the (classical) case
$\nu=1/2$ involving linear Brownian motion and the three-dimensional Bessel process
corresponds to the first-passage bridges used in \cite{Bet0} or \cite{BM}. We refrained from giving a more general
statement because our interest lies mainly in the case considered in the proposition. 

\proof Let us fix $0\leq \ve<x$. For every $t>0$ and $s\in[0,t)$, let $P^{s}_x$ be the law 
of $(R_u)_{0\leq u\leq s}$ under $\P^{(-1)}_x$. We define another probability measure 
$P^{\ve,s,t}_x$ on
$C([0,s],\R_+)$, which is absolutely continuous with respect to $P^{s}_x$, by letting the Radon-Nikodym derivative of $P^{\ve,s,t}_x$ with respect to $P^{s}_x$
be given by formula \eqref{RadonN} (note that $P^{\ve,s,t}_x$ is a probability measure by \eqref{Markov-densi}).
Then it is straightforward to verify that the collection $(P^{\ve,s,t}_x)_{t\in(s,\infty)}$ forms a regular version of the conditional distributions
of $(R_{u})_{0\leq u\leq s}$ knowing $T^{(R)}_\ve=t$, under $\P^{(-1)}_x(\cdot\cap\{T^{(R)}_\ve >s\})$. Furthermore, 
for every fixed $t>0$, the probability measures $P^{\ve,s,t}_x$ are consistent when $s$ varies, 
in the sense that, if $0\leq s<s'<t$, $P^{\ve,s,t}_x$ is the image of $P^{\ve,s',t}_x$ under the obvious
restriction mapping. It follows that we can define a
process $(X^{(t)}_u)_{0\leq u<t}$ with continuous sample paths on the time interval $[0,t)$ and such that for every $s\in[0,t)$,
the distribution of $(X^{(t)}_u)_{0\leq u\leq s}$ is $P^{\ve,s,t}_x$. From the Radon-Nikodym density \eqref{RadonN},
we can compute the finite-dimensional marginals of $X^{(t)}$,
\begin{align}
\label{fini-dim}
&\E[\varphi_1(X^{(t)}_{t_1})\varphi_2(X^{(t)}_{t_2})\cdots\varphi(X^{(t)}_{t_p})]\\
&=\frac{1}{r_t(x,\ve)}\int_{(\ve,\infty)^p} \dd y_1\ldots \dd y_p\,p^{(\ve)}_{t_1}(x,y_1)
p^{(\ve)}_{t_2-t_1}(y_1,y_2)\cdots p^{(\ve)}_{t_p-t_{p-1}}(y_{p-1},y_p)r_{t-t_p}(y_p,\ve)\,\varphi_1(y_1)\cdots\varphi_p(y_p),
\nonumber
\end{align}
for every $0<t_1<\cdots<t_p<t$ and every nonnegative measurable functions $\varphi_1,\ldots,\varphi_p$.
 We also set 
$X^{(t)}_t:=\ve$. Then it is not obvious that the sample paths of $X^{(t)}$ are continuous
at time $t$. To verify that this property holds, we use a time-reversal argument. It follows from
\eqref{Naga} and simple manipulations that, for every $0<t_1<\cdots<t_p$, the distributions of $(X^{(t)}_{t-t_1},X^{(t)}_{t-t_2},\ldots,X^{(t)}_{t-t_p})$
when $t$ varies in $(t_p,\infty)$ also form a regular version of the conditional distributions of 
$(R_{L^{(R)}_\ve+t_1},R_{L^{(R)}_\ve+t_2},\ldots,R_{L^{(R)}_\ve+t_p})$ knowing $L^{(R)}_x-L^{(R)}_\ve=t$, under $\P^{(5)}_0(\cdot\cap\{L^{(R)}_x-L^{(R)}_\ve>t_p\})$. 
By \eqref{fini-dim}, the density of $(X^{(t)}_{t-t_1},X^{(t)}_{t-t_2},\ldots,X^{(t)}_{t-t_p})$ is the function
$$(y_1,\ldots,y_p)\mapsto\frac{1}{r_t(x,\ve)} p^{(\ve)}_{t-t_p}(x,y_p)p^{(\ve)}_{t_p-t_{p-1}}(y_p,y_{p-1})\cdots p^{(\ve)}_{t_2-t_1}(y_2,y_1)
r_{t_1}(y_1,\ve).$$
If we integrate this density with respect to the measure $\mathbf{1}_{(t_p,\infty)}(t)r_t(x,\ve)\dd t$, we obtain that
the density of $(R_{L^{(R)}_\ve+t_1},R_{L^{(R)}_\ve+t_2},\ldots,R_{L^{(R)}_\ve+t_p})$ under $\P^{(5)}_0(\cdot\cap\{L^{(R)}_x-L^{(R)}_\ve>t_p\})$ is
$$(y_1,\ldots,y_p)\mapsto G^{(\ve)}(x,y_p)\,p^{(\ve)}_{t_p-t_{p-1}}(y_p,y_{p-1})\cdots p^{(\ve)}_{t_2-t_1}(y_2,y_1)
r_{t_1}(y_1,\ve).$$
By letting $x\to\infty$, we get that the density of $(R_{L^{(R)}_\ve+t_1},R_{L^{(R)}_\ve+t_2},\ldots,R_{L^{(R)}_\ve+t_p})$ under $\P^{(5)}_0$
is given by the same formula with $G^{(\ve)}(x,y_p)$ replaced by $G^{(\ve)}(\infty,y_p)$. 

From these finite-dimensional marginals distributions, we obtain that, for every $0<s<t$, the distribution 
of $(X^{(t)}_{t-u})_{0<u\leq s}$ is absolutely continuous with respect to the distribution 
of $(R_{L^{(R)}_\ve+u})_{0<u\leq s}$ under  $\P^{(5)}_0$, with a density given by
$$\w\mapsto \frac{p^{(\ve)}_{t-s}(x,\w(s))}{r_t(x,\ve)G^{(\ve)}(\infty,\w(s))}.$$
From this absolute continuity property, we get that $X^{(t)}_s\la \ve$ when $s\to t$, a.s.
So we can define a probability measure $\Pi^{(x,\ve)}_t(\dd \w)$ on $C([0,t],\R_+)$
as the distribution of $(X^{(t)}_s)_{0\leq s\leq t}$. It should be clear from our construction
that the collection $(\Pi^{(x,\ve)}_t)_{t>0}$ is a regular version of the
conditional distributions of $(R_s)_{0\leq s\leq T^{(R)}_\ve}$ knowing $T^{(R)}_\ve=t$ under $\P^{(-1)}_x$.

Finally, from the fact that $r_t(x+\ve,\ve)\la r_t(x,0)$ as $\ve\to 0$, and the analogous convergence $p^{(\ve)}_t(y,z)\la p_t(y,z)$,
which is derived from \eqref{density-condi}, it is a simple matter to verify that the finite-dimensional 
marginals of $\Pi^{(x+\ve,\ve)}_t$ converge to those of $\Pi^{(x,0)}_t$. Tightness of 
the collection $(\Pi^{(x+\ve,\ve)}_t)_{\ve\in(0,1)}$ is also easy
from the absolute continuity properties stated above. The last assertion of the proposition follows. \endproof

The probability measure $\Pi^{(x,0)}_t$ is also the law of 
the usual Bessel bridge of dimension $5$ from $x$ to $0$ over the time interval $[0,t]$. This may be 
verified from the finite-dimensional marginals in \eqref{fini-dim}, noting that the transition densities 
$p'_t(x,y)$, $t,x,y>0$, of the five-dimensional Bessel process satisfy $p'_t(x,y)=(\frac{y}{x})^3p_t(x,y)$ by the
the $h$-transform relation \eqref{h-transform} (we refer to \cite{PY}
for detailed information about Bessel bridges). 
The latter Bessel bridge
can be defined in a simpler way using the fact that the five-dimensional Bessel process
is the norm of a five-dimensional Brownian motion, and the additivity properties 
of squares of Bessel bridges (see e.g.~\cite{PY}). This interpretation also makes it possible 
to define a five-dimensional Bessel bridge $\mathbf{b}=(\mathbf{b}_t)_{0\leq t\leq 1}$ from $0$ 
to $0$ over the time interval $[0,1]$. The process $\mathbf{b}$ may indeed be obtained as the square root of the sum
of the squares of five independent standard (one-dimensional) Brownian bridges from $0$ to $0$
over $[0,1]$. One then easily verifies that the distribution of $\mathbf{b}_{1/2}$ has density
$$\rho(x)= \frac{64}{3\sqrt{\pi}}\,x^4\,e^{-2x^2}\;,\quad x>0.$$
Furthermore, conditionally on $\mathbf{b}_{1/2}=x$, the two processes 
$(\mathbf{b}_{\frac{1}{2}-t})_{0\leq t\leq 1/2}$ and $(\mathbf{b}_{\frac{1}{2}+t})_{0\leq t\leq 1/2}$
are independent and distributed according to $\Pi^{(x,0)}_{1/2}$ (that is, they are independent
Bessel bridges of dimension $5$ from $x$ to $0$).

We now turn to the main result of this section.

\begin{proposition}
\label{conv-bridge}
Let $\be=(\be_t)_{0\leq t\leq 1}$ be a normalized Brownian excursion. For every $\ve >0$, set 
$$C_\ve:=\E\Big[ 
\exp\Big(-\int_0^1 \frac{\dd t}{(\ve + \be_t)^2}\Big)\Big],$$
and write $\be^\ve=(\be^\ve_t)_{0\leq t\leq 1}$ for a random element of $C([0,1],\R_+)$ whose distribution is specified
by
$$\E[F(\be^\ve)]:= C_\ve^{-1}\,\E\Big[ F(\be)\,\exp\Big(-\int_0^1 \frac{\dd t}{(\ve + \be_t)^2}\Big)\Big],$$
for any nonnegative measurable function $F$ on  $C([0,1],\R_+)$.
Then, we have
\begin{equation}
\label{asympC}
\lim_{\ve\to 0} \ve^{-2}C_\ve=3,
\end{equation}
and
$$\be^\ve \build\la_{\ve\to 0}^{\rm(d)} \mathbf{b}$$
where $\mathbf{b}=(\mathbf{b}_t)_{0\leq t\leq 1}$ is a five-dimensional Bessel bridge from $0$ 
to $0$ over the time interval $[0,1]$. Finally, for every $\delta\in(0,\frac{1}{2})$, the total variation distance
between the distribution of $(\be^\ve_t+\ve)_{\delta\leq t\leq 1-\delta}$ and the 
distribution of $(\mathbf{b}_t)_{\delta\leq t\leq 1-\delta}$ tends to $0$ as $\ve\to 0$.
\end{proposition}

\rem As we will see later (cf.~formula \eqref{proba-ve} below), the quantity $C_\ve$ can also be interpreted as the probability in 
a free pointed Brownian disk of perimeter $1$ that the distance from the distinguished point to
the boundary is smaller than $\ve\sqrt{3}$. The asymptotics of $C_\ve$ when $\ve\to 0$ could 
therefore be derived from the distribution of the latter distance, which is known explicitly \cite{spine}.

\proof 
It is well known that, conditionally on $\be_{1/2}$, the two processes $(\be_{\frac{1}{2}+t})_{0\leq t\leq 1/2}$
and $(\be_{\frac{1}{2}-t})_{0\leq t\leq 1/2}$ are independent and follow the distribution of a linear Brownian
motion started from $\be_{1/2}$ and conditioned to hit $0$ for the first time at time $1/2$. This conditioned 
process can be defined in a way similar to Proposition \ref{Bessel-bridge}
(see Section 5.1 in \cite{Bet0} or Section 2.1 in \cite{BM}). Fix $x>0$ and $t>0$, and, for  every $s\in(0,t)$, 
write $\f_s$ for the $\sigma$-field on $C([0,t],\R)$ generated by $\w\mapsto(\w(r))_{0\leq r\leq s}$. Then
the Radon-Nikodym derivative on $\f_s$ of the law of 
Brownian motion started at $x$ and conditioned to hit $0$ for the first time at time $t$, with respect to
the law of Brownian motion started at $x$, is
$$\mathbf{1}_{\{\w(r)>0,\forall r\in[0,s]\}}\,\frac{q_{t-s}(\w(s),0)}{q_t(x,0)},$$
where the function $q_t(x,0)$ is as in \eqref{density-hittingBM} .

Write $(B_t)_{t\geq 0}$ for a linear Brownian motion
that starts at $x$ under the probability measure $\P_x$, and $T_y=\inf\{t\geq 0:B_t=y\}$
for every $y\in\R$. Let $F_1$ and $F_2$ be bounded Lipschitz continuous functions on $C([0,\frac{1}{2}],\R)$. We can summarize the 
first observation of the proof by the equality
\begin{align}
\label{conv-br-tech}
&\E\Big[ F_1\Big((\be_{\frac{1}{2}-t})_{0\leq t\leq \frac{1}{2}}\Big)F_2\Big((\be_{\frac{1}{2}+t})_{0\leq t\leq \frac{1}{2}}\Big)
\exp\Big(-\int_0^1 \frac{\dd t}{(\ve + \be_t)^2}\Big)\Big]\\
&= \int_0^\infty \dd x\,\pi(x)\,\E_x\Big[F_1\Big((B_t)_{0\leq t\leq \frac{1}{2}}\Big)\,e^{-\int_0^{1/2} \frac{\dd t}{(\ve+B_t)^2}}\,\Big|\, T_0=\frac{1}{2}\Big]
\E_x\Big[F_2\Big((B_t)_{0\leq t\leq \frac{1}{2}}\Big)\,e^{-\int_0^{1/2} \frac{\dd t}{(\ve+B_t)^2}}\,\Big|\, T_0=\frac{1}{2}\Big],\nonumber
\end{align}
where $\pi(x)= \frac{16}{\sqrt{\pi}}\,x^2\,e^{-2x^2}$ is the density of $\be_{1/2}$. We note that
\begin{equation}
\label{conv-br-tech2}
\E_x\Big[F_1\Big((B_t)_{0\leq t\leq \frac{1}{2}}\Big)\,e^{-\int_0^{1/2} \frac{\dd t}{(\ve+B_t)^2}}\,\Big|\, T_0=\frac{1}{2}\Big]
= \E_{x+\ve}\Big[F_1\Big((B_t-\ve)_{0\leq t\leq \frac{1}{2}}\Big)e^{-\int_0^{1/2} \frac{\dd t}{(B_t)^2}}\,\Big|\, T_\ve=\frac{1}{2}\Big].
\end{equation}
To study the right-hand side, we rely on the next lemma, where we use the notation introduced at the beginning of the section.

\begin{lemma}
\label{abso-cont}
Let $x>\ve>0$. Then, for every $t>0$ and every nonnegative measurable function $G$ on $C([0,t],\R)$,
\begin{equation}
\label{conv-br-tech4}
\E_x\Big[G((B_s)_{0\leq s\leq t})    e^{-\int_0^{t} \frac{\dd s}{(B_s)^2}}\,\Big|\,T_\ve=t\Big]
= \frac{\ve}{x} \frac{r_t(x,\ve)}{q_t(x,\ve)}\,\E^{(-1)}_x\Big[G((R_s)_{0\leq s\leq t}) \,\Big|\,T^{(R)}_\ve=t\Big].
\end{equation}
\end{lemma}

\proof We first prove that
\begin{equation}
\label{abso-tech}
\E_x\Big[ G((B_s)_{0\leq s\leq T_\ve})\,\exp\Big(-\int_0^{T_\ve} \frac{\dd s}{(B_t)^2}\Big)\Big]
= \ve^{-2}x^2\,\E^{(5)}_x\Big[\mathbf{1}_{\{T^{(R)}_\ve<\infty\}}\, G\Big((R_s)_{0\leq s\leq T^{(R)}_\ve}\Big)\Big],
\end{equation}
for every nonnegative measurable function $G$
on $\W$. 
This is basically a consequence of the absolute continuity relations between Bessel processes
(see e.g. \cite[Exercise XI.1.22]{RY}). 
These relations give the equality
$$\E_x\Big[ \mathbf{1}_{\{T_0>u\}}\,G((B_s)_{0\leq s\leq u})\,\exp\Big(-\int_0^{u} \frac{\dd s}{(B_s)^2}\Big)\Big]
=x^2\,\E^{(5)}_x\Big[(R_u)^{-2}\,G((R_s)_{0\leq s\leq u})\Big],$$
for every $u\geq 0$. So to get \eqref{abso-tech}, we just need to justify the replacement of the
constant time $u$ by the hitting time of $\ve$ in the last display. This can be done by standard approximation
techniques. For every $y>0$ and $n\geq 1$, write $[y]_n$ for the unique real of the form $k2^{-n}$, $k\in\N$,
such that $(k-1)2^{-n}<y\leq k2^{-n}$. 
Then, assuming that $G$ is bounded and continuous,
\begin{align*}
\E_x\Big[ G((B_s)_{0\leq s\leq T_\ve})\,e^{-\int_0^{T_\ve}  \frac{\dd s}{(B_s)^2}}\Big]
&=\lim_{n\to\infty} \E_x\Big[ 1_{\{[T_\ve]_n<T_0\}}\,G((B_s)_{0\leq s\leq [T_\ve]_n})\,e^{-\int_0^{[T_\ve]_n}  \frac{\dd s}{(B_s)^2}}\Big]\\
&=\lim_{n\to\infty} \sum_{k=1}^\infty \E_x\Big[ \mathbf{1}_{\{(k-1)2^{-n}<T_\ve\leq k2^{-n}<T_0\}}\,G((B_s)_{0\leq s\leq k2^{-n}})\,e^{-\int_0^{k2^{-n}}\! \frac{\dd s}{(B_s)^2}}\Big]\\
&=\lim_{n\to\infty} \sum_{k=1}^\infty x^2\,\E^{(5)}_x\Big[\mathbf{1}_{\{(k-1)2^{-n}<T^{(R)}_\ve\leq k2^{-n}\}}\, (R_{k2^{-n}})^{-2}\,G((R_s)_{0\leq s\leq k2^{-n}})
\Big]\\
&=\lim_{n\to\infty} x^2\,\E^{(5)}_x\Big[ (R_{[T^{(R)}_\ve]_n})^{-2}\,G\Big((R_s)_{0\leq s\leq [T^{(R)}_\ve]_n}\Big)
\,\mathbf{1}_{\{T^{(R)}_\ve<\infty\}}\Big]\\
&=\ve^{-2}\,x^2\,\E^{(5)}_x\Big[G\Big((R_s)_{0\leq s\leq T^{(R)}_\ve}\Big)\,\mathbf{1}_{\{T^{(R)}_\ve<\infty\}}\Big],
\end{align*}
where the use of dominated convergence in the last equality is justified by the fact that the
variable $(\inf_{s\geq 0}R_s)^{-2}$ is integrable under $\P^{(5)}_x$. This completes the proof of \eqref{abso-tech}.

Recalling that $\P^{(5)}_x(T^{(R)}_\ve<\infty)=(\ve/x)^3$ and using \eqref{h-trans}, we get that
the right-hand side of \eqref{abso-tech} can be written in the form
$$\frac{\ve}{x}\,\E^{(-1)}_x\Big[G\Big((R_s)_{0\leq s\leq T^{(R)}_\ve}\Big)\Big].$$

For $x>\ve\geq 0$, the density of $T_\ve$ under $\P_x$ is the function $t\mapsto q_t(x,\ve)$ defined in \eqref{density-hittingBM}. Also recall that the density of $T^{(R)}_\ve$
under $\P^{(-1)}_x$ is the function $t\mapsto r_t(x,\ve)$. It follows from \eqref{abso-tech} that, for $x>\ve>0$,
$$\int_0^\infty \dd t\,q_t(x,\ve)\,\E_x\Big[G((B_s)_{0\leq s\leq t})    e^{-\int_0^{t} \frac{\dd s}{(B_s)^2}}\,\Big|\,T_\ve=t\Big]
=\frac{\ve}{x}\int_0^\infty \dd t\,r_t(x,\ve)\,\E^{(-1)}_x\Big[G((R_s)_{0\leq s\leq t}) \,\Big|\,T^{(R)}_\ve=t\Big].  $$
Replacing $G((\w(s))_{0\leq s\leq t}) $ by $g(t)G((\w(s))_{0\leq s\leq t}) $, with an arbitrary nonnegative measurable
function $g$ on $\R_+$, we get that \eqref{conv-br-tech4} holds for 
Lebesgue almost every~$t>0$.  

To verify that \eqref{conv-br-tech4} indeed holds for {\it every} $t>0$, it is enough
to consider the special case where
$G((\w(s))_{s\leq t})=\mathbf{1}_{\{t>t_p\}}g_1(\w(t_1))\ldots g_p(\w(t_p))$ where $0<t_1<\cdots<t_p$ and $g_1,\ldots,g_p$
are bounded continuous functions from $\R$ into $\R_+$. Then formula \eqref{fini-dim} shows that the right-hand side
of \eqref{conv-br-tech4} is a continuous function of $t$ on $(t_p,\infty)$. On the other hand,
for $\delta>0$ and $t>t_p+\delta$, the 
left-hand side
of \eqref{conv-br-tech4} is bounded above by $I_\delta(t)$ and bounded below by $e^{-\delta/\ve^2}I_\delta(t)$,
where
\begin{align*}
I_\delta(t)&=\E_x\Big[g_1(B_{t_1})\ldots g_p(B_{t_p})\,e^{-\int_0^{t-\delta} \frac{\dd s}{(B_s)^2}}\,\Big|\,T_\ve=t\Big]\\
&=\E_x\Big[g_1(B_{t_1})\ldots g_p(B_{t_p})\,e^{-\int_0^{t-\delta} \frac{\dd s}{(B_s)^2}}\,\mathbf{1}_{\{T_\ve>t-\delta\}}\,
\frac{q_\delta(B_{t-\delta},\ve)}{q_t(x,\ve)}\Big].
\end{align*}
Thanks to dominated convergence, the last formula implies that $I_\delta(t)$ is a continuous function of $t\in(t_p+\delta,\infty)$. Letting $\delta\to 0$, it follows that the left-hand side
of \eqref{conv-br-tech4} is also a continuous function of $t$ on $(t_p,\infty)$. We conclude that 
\eqref{conv-br-tech4} holds for every $t>t_p$. This completes the proof of Lemma \ref{abso-cont}. \endproof

We return to the proof of Proposition \ref{conv-bridge}. 
From \eqref{conv-br-tech}, \eqref{conv-br-tech2} and \eqref{conv-br-tech4} (with $t=1/2$ and $x$ replaced by $x+\ve$), we get,
\begin{align}
\label{conv-br-tech8}
&\ve^{-2}\E\Big[ F_1\Big((\be_{\frac{1}{2}-t})_{0\leq t\leq 1/2}\Big)F_2\Big((\be_{\frac{1}{2}+t})_{0\leq t\leq 1/2}\Big)
\exp\Big(-\int_0^1 \frac{\dd t}{(\ve + \be_t)^2}\Big)\Big]\\
&=\!\int_0^\infty\! \dd x\,\frac{\pi(x)}{(x+\ve)^2}\frac{r_{\frac{1}{2}}(x+\ve,\ve)^2}{q_{\frac{1}{2}}(x+\ve,\ve)^2}
\E^{(-1)}_{x+\ve}\Big[F_1((R_s-\ve)_{0\leq s\leq \frac{1}{2}}) \Big|T^{(R)}_\ve\!=\!\frac{1}{2}\Big]\,
\E^{(-1)}_{x+\ve}\Big[F_2((R_s-\ve)_{0\leq s\leq \frac{1}{2}}) \Big|T^{(R)}_\ve\!=\!\frac{1}{2}\Big].\nonumber
\end{align}
We have the explicit expression
$$\frac{r_{\frac{1}{2}}(x+\ve,\ve)}{q_{\frac{1}{2}}(x+\ve,\ve)}= \frac{\sqrt{\pi}}{4\ve^3}\,\mathrm{erfc}(\frac{1}{2\ve}+x)
\,\exp\Big((\frac{1}{2\ve}+x)^2\Big) - \frac{1}{2\ve^2} + \frac{x}{\ve} +1,$$
from which it is a simple matter to get that
\begin{equation}
\label{conv-br-tech5}
\lim_{\ve\to 0} \frac{r_{\frac{1}{2}}(x+\ve,\ve)}{q_{\frac{1}{2}}(x+\ve,\ve)}
= \frac{r_{\frac{1}{2}}(x,0)}{q_{\frac{1}{2}}(x,0)}
= 2x^2,
\end{equation}
and 
 \begin{equation}
\label{conv-br-tech6}
\frac{r_{\frac{1}{2}}(x+\ve,\ve)}{q_{\frac{1}{2}}(x+\ve,\ve)} \leq K_1x^2 + K_2
\end{equation}
 with constants $K_1$ and $K_2$ that do not depend on $x>0$ and $\ve\in(0,1]$. Moreover, using the fact that
 $F_1$ and $F_2$ are Lipschitz continuous, we have for $i=1,2$, and for every $x>0$,
\begin{align}
\label{conv-br-tech7}
\lim_{\ve\to 0} \E^{(-1)}_{x+\ve}\Big[F_i((R_s-\ve)_{0\leq s\leq \frac{1}{2}}) \,\Big|\,T^{(R)}_\ve=\frac{1}{2}\Big]
 &=\lim_{\ve\to 0} \E^{(-1)}_{x+\ve}\Big[F_i((R_s)_{0\leq s\leq \frac{1}{2}}) \,\Big|\,T^{(R)}_\ve=\frac{1}{2}\Big]\nonumber\\
&=\E^{(-1)}_x\Big[F_i((R_s)_{0\leq s\leq \frac{1}{2}}) \,\Big|\,T^{(R)}_0=\frac{1}{2}\Big],
\end{align}
by the last assertion of Proposition \ref{Bessel-bridge}.

Thanks to \eqref{conv-br-tech5} and \eqref{conv-br-tech7}, we can now pass to the limit $\ve\to0$ in the 
right-hand side of \eqref{conv-br-tech8}, using \eqref{conv-br-tech6} to justify dominated convergence. It follows that
\begin{align}
\label{conv-br-tech3}
&\lim_{\ve\to 0}\ve^{-2}\E\Big[ F_1\Big((\be_{\frac{1}{2}-t})_{0\leq t\leq \frac{1}{2}}\Big)F_2\Big((\be_{\frac{1}{2}+t})_{0\leq t\leq \frac{1}{2}}\Big)
\exp\Big(-\int_0^1 \frac{\dd t}{(\ve + \be_t)^2}\Big)\Big]\nonumber\\
&\quad= 4\int_0^\infty \dd x\,\pi(x)\,x^2\,\E^{(-1)}_x\Big[F_1((R_s)_{0\leq s\leq \frac{1}{2}}) \,\Big|\,T^{(R)}_0=\frac{1}{2}\Big]
\,\E^{(-1)}_x\Big[F_2((R_s)_{0\leq s\leq \frac{1}{2}}) \,\Big|\,T^{(R)}_0=\frac{1}{2}\Big].
\end{align}
The particular case $F_1=F_2=1$ of \eqref{conv-br-tech3} gives
$$\lim_{\ve\to 0} \ve^{-2}C_\ve= 4 \int_0^\infty \dd x\,x^2\,\pi(x)= 3.$$
Furthermore, the function $x\mapsto \frac{4}{3}x^2\pi(x)$ is the density of $\mathbf{b}_{1/2}$, and it follows
from \eqref{conv-br-tech3} that
we have
\begin{align*}
&\lim_{\ve\to 0} (C_\ve)^{-1} \E\Big[ F_1\Big((\be_{\frac{1}{2}-t})_{0\leq t\leq \frac{1}{2}}\Big)F_2\Big((\be_{\frac{1}{2}+t})_{0\leq t\leq \frac{1}{2}}\Big)
\exp\Big(-\int_0^1 \frac{\dd t}{(\ve + \be_t)^2}\Big)\Big]\\
&\qquad= \E\Big[ F_1\Big((\mathbf{b}_{\frac{1}{2}-t})_{0\leq t\leq \frac{1}{2}}\Big)F_2\Big((\mathbf{b}_{\frac{1}{2}+t})_{0\leq t\leq \frac{1}{2}}\Big)\Big].
\end{align*}
This gives the convergence in distribution of $\be^\ve$ toward $\mathbf{b}$. 

It remains to prove the last assertion of the proposition. To this end, fix $\ve>0$, and let $A_1$ and $A_2$
be measurable subsets of $C([0,\frac{1}{2}-\delta],\R_+)$ such that $\min\{\w(t):0\leq t\leq \frac{1}{2}-\delta\}>\ve$
for every $\w\in A_1\cup A_2$.
Then, 
\begin{align*}
&\P\Big(\Big((\ve+\be^\ve_{\frac{1}{2}-t})_{0\leq t \leq \frac{1}{2}-\delta}, (\ve+\be^\ve_{\frac{1}{2}+t})_{0\leq t \leq \frac{1}{2}-\delta}\Big)\in A_1\times A_2\Big)\\
&\quad = C_\ve^{-1} \E\Big[ \mathbf{1}_{A_1}\big((\ve+\be_{\frac{1}{2}-t})_{0\leq t \leq \frac{1}{2}-\delta}\big)\,
 \mathbf{1}_{A_2}\big((\ve+\be_{\frac{1}{2}+t})_{0\leq t \leq \frac{1}{2}-\delta}\big)
 \exp\Big(-\int_0^1 \frac{\dd t}{(\ve + \be_t)^2}\Big)\Big]\\
 &\quad= C_\ve^{-1} \int_0^\infty \dd x\,\pi(x)\,\E_x\Big[\mathbf{1}_{A_1}\big((\ve+B_{t})_{0\leq t \leq \frac{1}{2}-\delta}\big)
  \exp\Big(-\int_0^{1/2} \frac{\dd t}{(\ve + B_t)^2}\Big)\,\Big|\,T_0=\frac{1}{2}\Big]\\
  &\hspace{35mm} \times
  \E_x\Big[\mathbf{1}_{A_2}\big((\ve+B_{t})_{0\leq t \leq \frac{1}{2}-\delta}\big)
  \exp\Big(-\int_0^{1/2} \frac{\dd t}{(\ve + B_t)^2}\Big)\,\Big|\,T_0=\frac{1}{2}\Big]\\
  &\quad = C_\ve^{-1} \int_\ve^\infty \dd x\,\pi(x-\ve)\,\E_x\Big[\mathbf{1}_{A_1}\big((B_{t})_{0\leq t \leq \frac{1}{2}-\delta}\big)
  \exp\Big(-\int_0^{1/2} \frac{\dd t}{(B_t)^2}\Big)\,\Big|\,T_\ve=\frac{1}{2}\Big]\\
  &\hspace{35mm} \times \E_x\Big[\mathbf{1}_{A_2}\big((B_{t})_{0\leq t \leq \frac{1}{2}-\delta}\big)
  \exp\Big(-\int_0^{1/2} \frac{\dd t}{(B_t)^2}\Big)\,\Big|\,T_\ve=\frac{1}{2}\Big].
\end{align*}
By \eqref{conv-br-tech4}, the quantities in the last display are equal to
\begin{align*}
&C_\ve^{-1}\ve^2\int_\ve^\infty \dd x\,\pi(x-\ve)\,\Big(\frac{r_{1/2}(x,\ve)}{x\,q_{1/2}(x,\ve)}\Big)^2\,
\P^{(-1)}_x\Big((R_{t})_{0\leq t \leq \frac{1}{2}-\delta}\in A_1\,\Big|\,T^{(R)}_{\ve}=\frac{1}{2}\Big)\\
&\hspace{60mm} \times\P^{(-1)}_x\Big((R_{t})_{0\leq t \leq \frac{1}{2}-\delta}\in A_2\,\Big|\,T^{(R)}_{\ve}=\frac{1}{2}\Big). 
\end{align*}
Recalling the Radon-Nikodym derivative \eqref{RadonN}, this is also equal to
$$C_\ve^{-1}\ve^2\int_\ve^\infty \dd x\,\frac{\pi(x-\ve)}{x^2q_{1/2}(x,\ve)^2}\,
\E^{(-1)}_x\Big[\mathbf{1}_{A_1}\big((R_{t})_{0\leq t \leq \frac{1}{2}-\delta}\big)r_\delta(R_{\frac{1}{2}-\delta},\ve)\Big]\,
\E^{(-1)}_x\Big[\mathbf{1}_{A_2}\big((R_{t})_{0\leq t \leq \frac{1}{2}-\delta}\big)r_\delta(R_{\frac{1}{2}-\delta},\ve)\Big].$$
It is convenient to consider that, under each probability measure $\P^{(-1)}_x$, we have an independent
copy $(R'_t)_{t\geq 0}$  of the Bessel process $(R_t)_{t\geq 0}$. The last display can then be written as 
$$C_\ve^{-1}\ve^2\int_\ve^\infty \dd x\,\frac{\pi(x-\ve)}{x^2q_{1/2}(x,\ve)^2}\,
\E^{(-1)}_x\Big[\mathbf{1}_{A_1}\big((R_{t})_{0\leq t \leq \frac{1}{2}-\delta}\big)
\mathbf{1}_{A_2}\big((R'_{t})_{0\leq t \leq \frac{1}{2}-\delta}\big)r_\delta(R_{\frac{1}{2}-\delta},\ve)
r_\delta(R'_{\frac{1}{2}-\delta},\ve)\Big].$$
We have thus obtained that
\begin{align}
\label{couple1}
&\P\Big(\Big((\ve+\be^\ve_{\frac{1}{2}-t})_{0\leq t \leq \frac{1}{2}-\delta}, 
(\ve+\be^\ve_{\frac{1}{2}+t})_{0\leq t \leq \frac{1}{2}-\delta}\Big)\in A\Big)\nonumber\\
&\quad =
C_\ve^{-1}\ve^2\int_\ve^\infty \dd x\,\frac{\pi(x-\ve)}{x^2q_{1/2}(x,\ve)^2}\,
\E^{(-1)}_x\Big[\mathbf{1}_{A}\big((R_{t})_{0\leq t \leq \frac{1}{2}-\delta},
(R'_{t})_{0\leq t \leq \frac{1}{2}-\delta}\big)\,r_\delta(R_{\frac{1}{2}-\delta},\ve)\,
r_\delta(R'_{\frac{1}{2}-\delta},\ve)\Big],
\end{align}
for any measurable subset $A$ of $C([0,\frac{1}{2}-\delta],(\ve,\infty))^2$.

On the other hand, by the last observation before the statement of the proposition, and using again
the Radon-Nidodym derivative \eqref{RadonN}, we have, for any measurable subset 
$A$ of $C([0,\frac{1}{2}-\delta],(0,\infty))^2$,
\begin{align}
\label{couple2}
&\P\Big(\Big((\mathbf{b}_{\frac{1}{2}-t})_{0\leq t \leq \frac{1}{2}-\delta}, 
(\mathbf{b}_{\frac{1}{2}+t})_{0\leq t \leq \frac{1}{2}-\delta}\Big)\in A\Big)\nonumber\\
&\quad = \frac{4}{3} \int_0^\infty \dd x\,\frac{x^2\pi(x)}{r_{1/2}(x,0)^2}\,
\E^{(-1)}_x\Big[\mathbf{1}_{A}\big((R_{t})_{0\leq t \leq \frac{1}{2}-\delta},
(R'_{t})_{0\leq t \leq \frac{1}{2}-\delta}\big)\,r_\delta(R_{\frac{1}{2}-\delta},0)\,
r_\delta(R'_{\frac{1}{2}-\delta},0)\Big].
\end{align}

By comparing the right-hand sides of \eqref{couple1} and \eqref{couple2}, we get that 
the total variation distance
between the distribution of $(\be^\ve_t+\ve)_{\delta\leq t\leq 1-\delta}$ and the 
distribution of $(\mathbf{b}_t)_{\delta\leq t\leq 1-\delta}$ is bounded above by
the sum of the quantities $\P(\min_{\delta\leq t\leq 1-\delta}\mathbf{b}_t \leq \ve) $ and
$$\int_\ve^\infty \dd x\,\E_x^{(-1)}\Bigg[\mathbf{1}_{\{m_{R,R'}>\ve\}}\Bigg|\frac{\ve^2\pi(x-\ve)}{C_\ve\,x^2q_{1/2}(x,\ve)^2}
r_\delta(R_{\frac{1}{2}-\delta},\ve)\,
r_\delta(R'_{\frac{1}{2}-\delta},\ve)- \frac{4x^2\pi(x)}{3r_{1/2}(x,0)^2}r_\delta(R_{\frac{1}{2}-\delta},0)\,
r_\delta(R'_{\frac{1}{2}-\delta},0)\Bigg|\Bigg]$$
where we have written $m_{R,R'}:=\min\{R_t\wedge R'_t:0\leq t\leq \frac{1}{2}-\delta\}$. 
Clearly, $\P(\min_{\delta\leq t\leq 1-\delta}\mathbf{b}_t \leq \ve) $ tends to $0$ as $\ve\to0$, 
so we need only check that the quantity in the last display also tends to $0$ as $\ve\to 0$.

We first note that
$$\frac{\ve^2\pi(x-\ve)}{C_\ve\,x^2q_{1/2}(x,\ve)^2}= \frac{\ve^2}{C_\ve}\,\frac{4\sqrt{\pi}}{x^2}\quad
\hbox{and}\quad \frac{4x^2\pi(x)}{3r_{1/2}(x,0)^2}=\frac{4\sqrt{\pi}}{3x^2}.$$
Recalling \eqref{asympC}, we see that the desired result will follow if we can prove that
\begin{equation}
\label{couple3}
\lim_{\ve \to 0} \int_\ve^\infty \frac{\dd x}{x^2}\,\E_x^{(-1)}\Big[\mathbf{1}_{\{m_{R,R'}>\ve\}}\,\Big|r_\delta(R_{\frac{1}{2}-\delta},\ve)\,
r_\delta(R'_{\frac{1}{2}-\delta},\ve)- r_\delta(R_{\frac{1}{2}-\delta},0)\,
r_\delta(R'_{\frac{1}{2}-\delta},0)\Big|\Big] =0,
\end{equation}
and
\begin{equation}
\label{couple4}
\int_0^\infty \frac{\dd x}{x^2}\,\E_x^{(-1)}\Big[ \mathbf{1}_{\{m_{R,R'}>0\}}r_\delta(R_{\frac{1}{2}-\delta},0)\,
r_\delta(R'_{\frac{1}{2}-\delta},0)\Big]<\infty.
\end{equation}
The proof of \eqref{couple4} is immediate since the integral in \eqref{couple4} is equal, up to a multiplicative 
constant, to the right-hand side of \eqref{couple2} with $A=C([0,\frac{1}{2}-\delta],(0,\infty))^2$.
Then, we observe that the integral in \eqref{couple3} is equal to
$$ \int_\ve^\infty \frac{\dd x}{x^2}\,\int_\ve^\infty\dd y\int_\ve^\infty \dd z \,p^{(\ve)}_{\frac{1}{2}-\delta}(x,y)p^{(\ve)}_{\frac{1}{2}-\delta}(x,z)\,\Big|r_\delta(y,\ve)\,
r_\delta(z,\ve)- r_\delta(y,0)\,
r_\delta(z,0)\Big|.$$
and we can bound $p^{(\ve)}_{\frac{1}{2}-\delta}(x,y)p^{(\ve)}_{\frac{1}{2}-\delta}(x,z)$ by
$p_{\frac{1}{2}-\delta}(x,y)p_{\frac{1}{2}-\delta}(x,z)$. Furthermore, we know that, for every fixed $y,z>0$, the 
quantities $|r_\delta(y,\ve)\,
r_\delta(z,\ve)- r_\delta(y,0)\,
r_\delta(z,0)|$ tend to $0$ as $\ve\to 0$, and these quantities (for $y>\ve$ and $z>\ve$) are uniformly bounded by
a constant depending only on $\delta$ (this follows from our explicit formula for $r_t(x,\ve)$). 
In order to justify dominated convergence and to get \eqref{couple3}, it suffices to verify that
\begin{equation}
\label{couple5}
\int_0^\infty \frac{\dd x}{x^2}\,\int_0^\infty\dd y\int_0^\infty \dd z \,p_{\frac{1}{2}-\delta}(x,y)p_{\frac{1}{2}-\delta}(x,z)<\infty.
\end{equation}
However, using \eqref{h-trans} and writing $p'_t(x,y)$ for the transition densities of
the five-dimensional Bessel process, we have for every $x>0$,
$$\int_0^\infty\dd y\,p_{\frac{1}{2}-\delta}(x,y)= \int_0^\infty\dd y\,\frac{x^3}{y^3}\,p'_{\frac{1}{2}-\delta}(x,y)
= x^3\,\E_x^{(5)}[(R_{\frac{1}{2}-\delta})^{-3}]\leq x^3\,(x^{-3}\wedge K),$$
with a constant $K$ depending only on $\delta$. It follows that the integral in \eqref{couple5} is bounded
above by $\int_0^\infty \dd x\,x^{-2}\,(1\wedge Kx^{3})^2<\infty$. 
This completes the proof of \eqref{couple3} and of the proposition.
\endproof

\section{The free pointed Brownian disk and the Brownian half-plane}
\label{sec:BD-BHP}

In this section we recall the definitions of the (free pointed) Brownian disk and of the Brownian half-plane
along the lines of \cite{BMR,Bet,BM,GM0}. Our presentation is a little different from the latter papers
and better suited to our applications.

\medskip\noindent
{\bf The free pointed Brownian disk.} We consider a Poisson point measure $\nn=\sum_{i\in I} \delta_{(t_i,\omega_i)}$ on
$[0,1]\times \S$ with intensity
$$2\,\mathbf{1}_{[0,1]}(t)\,\dd t\,\N_0(\dd \omega).$$
We then introduce the compact metric space $\mathfrak{H}$, which is obtained from
the disjoint union
\begin{equation}
\label{tree-disk}
[0,1] \cup \Big(\bigcup_{i\in I} \t_{(\omega_i)}\Big)
\end{equation}
by identifying $0$ with $1$ and, for every $i\in I$, the root $\rho_{(\omega_i)}$ of $\t_{(\omega_i)}$
with the point $t_i$ of $[0,1]$. The metric $\dd_\mathfrak{H}$ on $\mathfrak{H}$ is defined as follows. 
First, the restriction of $\dd_\mathfrak{H}$ to each tree $\t_{(\omega_i)}$ is 
the metric $d_{(\omega_i)}$. Then, if $u,v\in[0,1]$, we take $\dd_\mathfrak{H}(u,v)=\min\{u\vee v-u\wedge v, 1-u\vee v+u\wedge v)$. If $u\in [0,1]$, and $v\in\t_{(\omega_i)}$ for some $i\in I$, $\dd_\mathfrak{H}(u,v)=\dd_\mathfrak{H}(u,\rho_{(\omega_i)})+ d_{(\omega_i)}(\rho_{(\omega_i)},v)$. Finally
  if $u\in \t_{(\omega_i)}$ and $v\in \t_{(\omega_j)}$, with $j\not =i$, 
$$\dd_\mathfrak{H}(u,v)=d_{(\omega_i)}(u,\rho_{(\omega_i)})+d_\mathfrak{H}(\rho_{(\omega_i)},\rho_{(\omega_j)})+ d_{(\omega_j)}(\rho_{(\omega_j)},v).$$
The volume measure on $\mathfrak{H}$ is just the sum of the volume measures on the trees $\t_{(\omega_i)}$, $i\in I$.

If $\Sigma:=\sum_{i\in I}\sigma(\omega_i)$ is the total mass of the volume measure, we can define a cyclic clockwise exploration $(\ee_t)_{0\leq t\leq\Sigma}$
of $\mathfrak{H}$, informally by concatenating the mappings $p_{(\omega_i)}:[0,\sigma(\omega_i)]\la \t_{(\omega_i)}$ in the 
order prescribed by the $t_i$'s. To give a more precise definition, set
$$\beta_s:= \sum_{i\in I} \mathbf{1}_{\{t_i\leq s\}}\,\sigma(\omega_i)\,,\ \beta_{s-}:= \sum_{i\in I} \mathbf{1}_{\{t_i< s\}}\,\sigma(\omega_i)\,,$$
for every $s\in[0,1]$. 
Then, for every $t\in[0,\Sigma]$, we define $\ee_t\in\mathfrak{H}$ as follows. We observe that there is a unique $s\in[0,1]$ such that $\beta_{s-}\leq t\leq \beta_{s}$, and:
\begin{description}
\item[$\bullet$] Either there is a (unique) $i\in I$ such that $s=t_i$, and we set
$\ee_t:= p_{(\omega_i)}(t-\beta_{t_i-})$.
\item[$\bullet$] Or there is no such $i$ and we set $\ee_t:=s$.
\end{description}
Note that $\ee_\Sigma=\ee_0$ (because $1$ is identified to $0$ in $\mathfrak{H}$). 

The clockwise exploration allows us to define ``intervals''
in $\mathfrak{H}$. Let us make the convention that, if 
$s,t\in[0,\Sigma]$ and $s>t$, the (real) interval $[s,t]$ is defined by $[s,t]:=[s,\Sigma]\cup [0,t]$ (of course, if 
$s\leq t$, $[s,t]$ is the usual interval).  
Then, for every $u,v\in\mathfrak{H}$, such that $u\not =v$, there is a smallest interval $[s,t]$, with $s,t\in[0,\Sigma]$, such that
$\ee_s=u$ and $\ee_t=v$, and we define  
$$[|u,v|]:=\{\ee_r:r\in[s,t]\}.$$
We have typically $[|u,v|]\not =[|v,u|]$. Of course, we take $[|u,u|]=\{u\}$. Note that we use the notation 
$[|u,v|]$ rather than $[u,v]$ to avoid confusion with intervals of the real line.

We next assign real labels to the points of $\mathfrak{H}$. To this end, we let $(\be_t)_{0\leq t\leq 1}$ be a normalized
Brownian excursion, which is independent of $\nn$. For $t\in[0,1]$,
we set $\Lambda_t:=\sqrt{3}\,\be_t$, and for $u\in\t_{(\omega_i)}$, $i\in I$, 
$$\Lambda_u:=\sqrt{3}\,\be_{t_i} + \ell_u(\omega_i),$$
where we recall that $\ell_u(\omega_i)$ is the label of $u$ in $\t_{(\omega_i)}$. 
By \cite[Lemma 11]{Bet}, $\min\{\Lambda_u:u\in \mathfrak{H}\}$ is attained at a unique point $v_*$ of $\mathfrak{H}$,
and we set $\Lambda_*:=\Lambda_{v_*}$ to simplify notation.

Labels allow us to define the pseudo-metric $D$ on $\mathfrak{H}$ as follows. For every $u,v\in\mathfrak{H}$,
we first set
\begin{equation}
\label{pseudo-BD1}
D^{\circ}(u,v):=\Lambda_u + \Lambda_v -2 \max\Big( \inf_{w\in[|u,v|]} \Lambda_w,\inf_{w\in[|v,u|]} \Lambda_w\Big),
\end{equation}
and then
\begin{equation}
\label{pseudo-BD2}
D(u,v) := \inf_{u_0=u,u_1,\ldots,u_p=v} \sum_{i=1}^p D^{\circ}(u_{i-1},u_i),
\end{equation}
where the infimum is over all choices of the integer $p\geq 1$ and of the
finite sequence $u_0,u_1,\ldots,u_p$ in $\mathfrak{H}$ such that $u_0=u$ and
$u_p=v$. One immediately verifies that $D(u,v)\geq |\Lambda_u-\Lambda_v|$ for every $u,v\in\mathfrak{H}$. 
It easily follows that, for every $u\in\mathfrak{H}$, $D(u,v_*)=D^\circ(u,v_*)= \Lambda_u-\Lambda_{*}$.
We also notice that the mapping $(u,v)\mapsto D(u,v)$ is continuous on $\mathfrak{H}\times\mathfrak{H}$
(note that $D^{\circ}(u_n,u)\la 0$ if $u_n\to u$ in $\mathfrak{H}$, and use the triangle inequality).

We abuse notation by writing $\mathfrak{H}/\{D=0\}$ for the quotient space of $\mathfrak{H}$ with respect 
to the equivalence relation 
defined by setting $u\sim v$ if and only if $D(u,v)=0$.  

\begin{definition}
\label{def-disk}
The free pointed Brownian disk with perimeter $1$ is the quotient space $\D^\bullet:=\mathfrak{H}/\{D=0\}$, which is equipped
with the distance induced by $D$ and with a distinguished point which is the equivalence class of $v_*$. The volume
measure on $\D^\bullet$ is the pushforward of the volume measure on $\mathfrak{H}$
under the canonical projection.
\end{definition}

We may and will view $\D^\bullet$ as a (random) pointed compact measure metric space, that is, as an element
of the space $\M^\bullet$ of Section \ref{sec:GHP}. The reader will easily check that this presentation of the free pointed
Brownian disk is consistent with the one in \cite{BM}. Note that the role of the Brownian excursion $\be$
is played in \cite{Bet,BM} by a standard Brownian bridge. The celebrated Vervaat transformation \cite{Vervaat} connecting the
Brownian bridge with the Brownian excursion shows that this makes no difference (note that adding a
random constant to all labels does not change the definition of $D$). 

We will use the notation $\Pi$ for the canonical projection
from $\mathfrak{H}$ onto $\D^\bullet$. We note that $\D^\bullet$ is a length space. This can be verified by
observing that, for
every $u,v\in\mathfrak{H}$, $D^\circ(u,v)$ is the length of a continuous curve from $\Pi(u)$
to $\Pi(v)$ in $\D$, namely the curve obtained by concatenating the respective simple geodesics 
from $\Pi(u)$ and $\Pi(v)$ to $\Pi(v_*)$ until the point where they merge (see e.g.~\cite[Section 2.6]{Uniqueness} for a definition of simple geodesics in the
Brownian map, which is immediately adapted to the present setting).

Furthermore, the space $\D^\bullet$ is homeomorphic to the closed unit disk of the plane \cite{Bet}, and, in this homeomorphism,
the unit circle corresponds to $\partial\D^\bullet:=\Pi([0,1))$. We define the uniform measure $\mu$ on $\partial\D^\bullet$
as the pushforward of Lebesgue measure on $[0,1)$ under $\Pi$. From \cite[Lemma 14]{Bet}, 
one knows that a.s. for every $u\in [0,1)\subset\mathfrak{H}$, the equivalence class of $u$ in the 
quotient $\mathfrak{H}/\{D=0\}$ is a singleton (no point of $[0,1)$ is identified with another point of $\mathfrak{H}$).
Notice that \cite{Bet} deals with the slightly different model where the total volume of $\mathfrak{H}$ is fixed
(corresponding to the Brownian disk with fixed volume and perimeter), but
the result also applies to our setting. In particular the mapping $[0,1)\ni u\mapsto \Pi(u)$ is injective. 

We will keep the notation $D$ for the metric of $\D^\bullet$ and, without risk of confusion, we identify $v_*$ with
$\Pi(v_*)$. For every $x\in\D^\bullet$, we set $\Lambda_x:=\Lambda_u$, where $u$ is
a point of $\mathfrak{H}$ such that $\Pi(u)=x$ (the bound $|\Lambda_u-\Lambda_v|\leq D(u,v)$ shows that this does 
not depend on the choice of $u$). Then, we have $D(v_*,x)=\Lambda_x-\Lambda_*$ for every $x\in\D^\bullet$.

The fact that the Brownian bridge used in the presentation of \cite{Bet,BM} is replaced here by
a Brownian excursion has an important consequence. In \cite{Bet,BM}, the 
equivalence class of $0$ is a typical point of the boundary, in a sense
that can be made precise,
whereas here $\Pi(0)$ is the point of $\partial\D^\bullet$ that is closest to $v_*$ (and is therefore a very special point). 
The distance from the distinguished point $v_*$ to $\partial\D^\bullet$ is
$$D(v_*,\partial\D^\bullet)=
D(v_*,\Pi(0))=-\Lambda_* .$$
The explicit distribution of $D(v_*,\partial\D^\bullet)$ is given in \cite{spine}, but for us it will be sufficient to know
the asymptotics of $\P(D(v_*,\partial\D^\bullet)\leq\ve)$ when $\ve>0$. To this end, 
note that the event $\{D(v_*,\partial\D^\bullet)\leq\ve\}=\{\Lambda_*\geq-\ve\}$ occurs if and only if we have $\sqrt{3}\be_{t_i}+W_*(\omega_i)\geq -\ve$
for every atom $(t_i,\omega_i)$ of $\nn$. Using \eqref{hittingpro}, we obtain
that 
\begin{equation}
\label{proba-ve}
\P(D(v_*,\partial\D^\bullet)\leq\ve)= \E\Big[\exp\Big(-3\int_0^1 \frac{\dd t}{(\sqrt{3}\,\be_t+ \ve)^2}\Big)\Big]
= \E\Big[\exp\Big(-\int_0^1 \frac{\dd t}{(\be_t+ \ve/\sqrt{3})^2}\Big)\Big]
\end{equation}
and \eqref{asympC} then yields
\begin{equation}
\label{asymp-pro}
\lim_{\ve \to 0} \ve^{-2}\,\P(D(v_*,\partial\D^\bullet)\leq\ve) = 1.
\end{equation}

\smallskip

\noindent{\bf The Brownian half-plane.} 

We now present the construction of the Brownian half-plane along the 
lines of \cite{BMR} or \cite{GM0}. This construction is very similar to that
of the free pointed Brownian disk presented above, and
we will therefore omit a few details.

We consider a Poisson point measure $\nn_\infty=\sum_{j\in J} \delta_{(t^\infty_j,\omega^\infty_j)}$ on
$\R\times \S$ with intensity $2\,\dd t\,\N_0(\dd \omega)$. 
We introduce the locally compact metric space $\mathfrak{H}_\infty$, which is obtained from
the disjoint union
\begin{equation}
\label{tree-disk-infty}
\R \cup \Big(\bigcup_{j\in J} \t_{(\omega^\infty_j)}\Big)
\end{equation}
by identifying, for every $j\in J$, the root $\rho_{(\omega^\infty_j)}$ of $\t_{(\omega_j)}$
with the point $t^\infty_j$ of $\R$. The metric $\dd_{\mathfrak{H}_\infty}$ on $\mathfrak{H}_\infty$ is defined 
in the same way as the metric $\dd_{\mathfrak{H}}$ on $\mathfrak{H}$ was defined above
(the restriction of $\dd_{\mathfrak{H}_\infty}$ to $\R$ is the
usual Euclidean metric). We note that $\mathfrak{H}_\infty$ is a (non-compact) $\R$-tree.
The volume measure on $\mathfrak{H}_\infty$ is the $\sigma$-finite measure that puts no mass on $\R$
and whose restriction to each tree $\t_{(\omega_j)}$, $j\in J$, is the volume measure on this tree.

Similarly as above, the clockwise exploration $(\ee^\infty_s)_{s\in\R}$ 
of $\mathfrak{H}_\infty$ is defined by concatenating the mappings $p_{(\omega^\infty_j)}:[0,\sigma(\omega^\infty_j)]\la \t_{(\omega^\infty_j)}$ in the 
order prescribed by the $t^\infty_j$'s, in such a way that $\ee^\infty_0=0$ (so the points $(\ee^\infty_s,s\geq 0)$ correspond
exactly to the union of $\R_+$ and of the trees $\t_{(\omega^\infty_j)}$ for indices $j$ such that $t^\infty_j\geq 0$). 
We omit the precise description of $(\ee^\infty_s)_{s\in\R}$, which should be obvious from the analogous
definition of $(\ee_s)_{s\in[0,\Sigma]}$ given above.

The clockwise exploration allows us to define intervals 
in the space $\mathfrak{H}_\infty$. We now make the convention that, if 
$s,t\in\R$ and $s>t$, $[s,t]=[s,\infty)\cup(-\infty,t]$.  
Then, for every $u,v\in\mathfrak{H}_\infty$, such that $u\not =v$, 
we set $[|u,v|]:=\{\ee^\infty_r:r\in[s,t]\}$, where $[s,t]$ is the smallest  ``interval'' such that
$\ee^\infty_s=u$ and $\ee^\infty_t=v$. Note that at least one
of the two intervals $[|u,v|]$ and $[|v,u|]$ is compact.

In order to assign real labels to the points of $\mathfrak{H}_\infty$, we consider a two-sided Brownian motion
$B=(B_t)_{t\in\R}$ (in other words, $(B_{t})_{t\geq0}$ and $(B_{-t})_{t\geq 0}$
are two independent linear Brownian motions started from $0$), and 
we assume that $B$ is independent of $\nn_\infty$. We
set $\Lambda^\infty_u:=\sqrt{3}\,B_u$ if $u\in \R$, and 
for $u\in\t_{(\omega^\infty_j)}$, $j\in J$, 
$$\Lambda^\infty_u:=\sqrt{3}\,B_{t^\infty_j} + \ell_u(\omega^\infty_j).$$
Then, for every $u,v\in\mathfrak{H}_\infty$,
we set
\begin{equation}
\label{pseudo-infty1}
D^{\infty,\circ}(u,v):=\Lambda^\infty_u + \Lambda^\infty_v -2 \max\Big( \inf_{w\in[|u,v|]} \Lambda^\infty_w,\inf_{w\in[|v,u|]} \Lambda^\infty_w\Big),
\end{equation}
and 
\begin{equation}
\label{pseudo-infty2}
D^\infty(u,v) := \inf_{u_0=u,u_1,\ldots,u_p=v} \sum_{i=1}^p D^{\infty,\circ}(u_{i-1},u_i)
\end{equation}
where the infimum is over all choices of the integer $p\geq 1$ and of the
finite sequence $u_0,u_1,\ldots,u_p$ in $\mathfrak{H}_\infty$ such that $u_0=u$ and
$u_p=v$. It is immediate that
$D^\infty(u,v)\geq |\Lambda^\infty_u-\Lambda^\infty_v|$ for every $u,v\in\mathfrak{H}_\infty$, and
we also note that the mapping $(u,v)\mapsto D^\infty(u,v)$ is continuous on $\mathfrak{H}_\infty\times\mathfrak{H}_\infty$.
Furthermore, the so-called cactus bound states that, for every $u,v\in\mathfrak{H}_\infty$, 
\begin{equation}
\label{cactus}
D^\infty(u,v)\geq \Lambda^\infty_u + \Lambda^\infty_v- 2\min_{w\in\llbracket u,v\rrbracket}\Lambda^\infty_w,
\end{equation}
where $\llbracket u,v\rrbracket$ denotes the geodesic segment between $u$ and $v$ in the $\R$-tree $\mathfrak{H}_\infty$
(not to be confused with the interval $[|u,v|]$). See formula (4) in \cite{Plane} for a short proof
in the case of the Brownian map, which is immediately extended to the present setting.

We finally notice that, if $u,v\in\mathfrak{H}_\infty$ are such that $[|u,v|]=\{\ee^\infty_r:r\in[s,t]\}$
with $s>t$, then trivially $\inf_{w\in[|u,v|]} \Lambda^\infty_w=-\infty$ and thus the maximum 
in \eqref{pseudo-infty1} is equal to $\inf_{w\in[|v,u|]} \Lambda^\infty_w$ (this occurs in
particular if $u\in \t_{(\omega_i)}$ and $v\in \t_{(\omega_j)}$ with $t_j<t_i$). 

\begin{definition}
\label{def-half-plane}
The Brownian half-plane is the quotient space $\H:=\mathfrak{H}_\infty/\{D^\infty=0\}$, which is equipped
with the distance induced by $D^\infty$ and with a distinguished point which is the equivalence class of $0$. The volume
measure on $\H$ is the pushforward of the volume measure on $\mathfrak{H}_\infty$
under the canonical projection.
\end{definition}

We use the notation $\Pi_\infty$ for the canonical projection
from $\mathfrak{H}_\infty$ onto $\H$, and keep the notation $D^\infty$ for the metric on $\H$. By the same argument as for the Brownian disk, $\H$ is a length space,
and (for instance by using Lemma \ref{dist-infty} below) it is easy to verify that closed balls 
in $\H$ are compact. Thus we may and will view $\H$ as a random element of the space $\M^\bullet_{bcl}$ of Section \ref{sec:GHP}.

The space $\H$ is homeomorphic to the usual upper half-plane  and in this homeomorphism
the real line corresponds to $\partial\H:=\Pi_\infty(\R)$ (see \cite[Section 1.5]{GM0}). The uniform measure $\mu_\infty$ on $\partial\H$
is defined as the pushforward of Lebesgue measure on $\R$ under $\Pi_\infty$. For our purposes, it is important to
note that the equivalence class of any $u\in\R$ in the quotient $\H=\mathfrak{H}_\infty/\{D^\infty=0\}$
is a singleton (no point of $\R$ is identified to another point of $\mathfrak{H}_\infty$). Indeed, comparing the constructions of the Brownian half-plane and of the free pointed Brownian disk
(and also using Lemma \ref{dist-infty} below), one may observe that the existence of a pair $\{u,v\}$ of distinct points of $\mathfrak{H}_\infty$ such that $u\in\R$ and $D^\infty(u,v)=0$ would imply the existence of a pair with similar properties in the Brownian disk, and we know from
\cite{Bet} that this does not occur. 

For reals $a<b$, we set
$$\mathfrak{H}_\infty^{[a,b]}:=[a,b]\cup \Big(\bigcup_{j\in J, a\leq t_j\leq b} \t_{(\omega_j)}\Big)$$
of course with the same identifications as in the definition of $\mathfrak{H}_\infty$, so that
$\mathfrak{H}_\infty^{[a,b]}$ is a subset of $\mathfrak{H}_\infty$.

\begin{lemma}
\label{dist-infty}
We have a.s.
$$\lim_{a\to\infty} \Big(\inf_{u\in\mathfrak{H}_\infty\backslash \mathfrak{H}_\infty^{[-a,a]}} D^\infty(0,u)\Big)=\infty.$$
\end{lemma}

\proof Let $u\in \mathfrak{H}_\infty\backslash \mathfrak{H}_\infty^{[-a,a]}$, and assume for definiteness that 
$u$ belongs to a tree $\t_{(\omega_j)}$ with $t_j>a$, or that $u\in(a,\infty)$. The cactus bound \eqref{cactus}
ensures that 
$$D^\infty(0,u) \geq \Lambda^\infty_u- 2\min_{v\in\llbracket 0,u\rrbracket}\Lambda^\infty_v
\geq -\min_{v\in\llbracket 0,u\rrbracket}\Lambda^\infty_v.$$
Since the geodesic segment $\llbracket 0,u\rrbracket$
contains the interval $[0,a]$, the right-hand side of the preceding display is bounded
below by $-\min_{0\leq t\leq a}B_t$. Finally, the infimum in the lemma is bounded below by
$$\Big(-\min_{0\leq t\leq a}B_t\Big)\wedge\Big(-\min_{-a\leq t\leq 0}B_t\Big)$$
which tends to $\infty$ as $a\to\infty$. \endproof

\noindent{\bf Notation.} Without risk of confusion, we will use the same notation $\mathrm{Vol}(\cdot)$ for the
volume measure on any of the spaces $\mathfrak{H},\mathfrak{H}_\infty,\D^\bullet$ and $\H$, as well as on the space
$\D$ introduced below. 

\section{Approximating the uniform measure on the boundary}
\label{sec:approx}

We consider the free pointed Brownian disk $(\D^\bullet,D)$ as defined in the previous section. Recall the
definition of the uniform probability measure $\mu$ on $\partial \D^\bullet$ as the pushforward
of Lebesgue measure on $[0,1)$ under the canonical projection $\Pi$.
The goal of this section is to prove the following useful approximation result.

\begin{theorem}
\label{approx-unif-bdry}
For every $\ve>0$, let $\mu_\ve$ be the finite measure on $\D^\bullet$ defined by
$$\langle \mu_\ve,\varphi\rangle = \ve^{-2}\,\int_{\D^\bullet} \mathrm{Vol}(\dd x)\,
\mathbf{1}_{\{D(x,\partial\D^\bullet)\leq\ve\}}\,\varphi(x).$$
Then a.s. $\mu_\ve$ converges weakly to $\mu$ as $\ve\to 0$.
\end{theorem}

It is proved in \cite[Proposition 2]{Disks} that the measures $\mu_\ve$ converge weakly to
a probability measure $\nu$, which is also called the uniform probability measure on the boundary in \cite{Disks}.
So Theorem \ref{approx-unif-bdry} is equivalent to the statement $\mu=\nu$. Unfortunately, this
equality is not easy to prove, because the construction of the free pointed Brownian disk in
\cite{Disks} is very different from the one presented in Section \ref{sec:BD-BHP} (see the
comments in the introduction of \cite{Disks}). So below, we will essentially prove 
the convergence of $\mu_\ve$ to $\mu$ independently of the results of \cite{Disks} --- we
still need these results to get the value of the constant $\kappa$ that appears in Lemma \ref{voisin-infty}
below. 

Before we proceed to the proof of Theorem \ref{approx-unif-bdry}, we need a few preliminary lemmas, 
which are mainly concerned with the case of
the Brownian half-plane $(\H,D^\infty)$ constructed in the previous section as a quotient space of $\mathfrak{H}_\infty$.
As
previously, we view $\R$ as a subset of $\mathfrak{H}_\infty$. Recall the notation $\mathfrak{H}_\infty^{[a,b]}$
introduced before Lemma \ref{dist-infty}.

\begin{lemma}
\label{voisin-infty}
There exists a constant $\kappa\in(0,\infty]$ such that, for any 
reals $a,b$ with $a<b$, we have
$$\ve^{-2}\,\mathrm{Vol}(\{u\in \mathfrak{H}_\infty^{[a,b]}:D^\infty(u,\R)\leq \ve\}) \build{\la}_{\ve\to 0}^{} \kappa(b-a),$$
in probability.
\end{lemma}

\rem We will see later that $\kappa=1$, but at the present stage, we do not exclude the 
possibility that $\kappa=\infty$. 

\proof Consider first the case $a=0,b=1$. Simple arguments relying on the invariance of $B$ and $\N_0$
under scaling transformations show that
$$\mathrm{Vol}(\{u\in \mathfrak{H}_\infty^{[0,1]}:D^\infty(u,\R)\leq \ve\}) 
\build{=}_{}^{\rm(d)} \ve^4\,\mathrm{Vol}(\{u\in \mathfrak{H}_\infty^{[0,1/\ve^2]}:D^\infty(u,\R)\leq 1\}).$$
So we will get the desired convergence for $a=0,b=1$ if we can verify that
\begin{equation}
\label{voisin-tech}
\frac{1}{n} \,\mathrm{Vol}(\{u\in \mathfrak{H}_\infty^{[0,n]}:D^\infty(u,\R)\leq 1\}) \build{\la}_{\ve\to 0}^{} \kappa,
\end{equation}
in probability, with some constant $\kappa\in(0,\infty]$. To this end, we may assume that
the pair $(B,\nn_\infty)$ is defined on the canonical space $\Omega_\circ:=C(\R,\R)\times M_p(\R\times \S)$,
in such a way that $B_t(\w,\gamma)=\w(t)$ and $\nn_\infty(\w,\gamma)=\gamma$
for $(\w,\gamma)\in \Omega_\circ$. The space $\Omega_\circ$ is equipped with the
unique probability measure $\P$ under which $B$ and $\nn_\infty$
have the required properties. The shift $\theta$ on $\Omega_\circ$ is then defined by
$$\theta\Big(\w,\sum_{k\in I}\delta_{(t_k,\omega_k)}\Big) = \Big(\w(1+\cdot)-\w(1), \sum_{k\in I}\delta_{(t_k-1,\omega_k)}\Big),$$
and $\P$ is invariant under $\theta$. 

For every integers $i<j$, set 
$$V_{i,j}:=\mathrm{Vol}(\{u\in \mathfrak{H}_\infty^{[i,j]}:D^\infty(u,\R)\leq 1\}).$$
Then,
$$V_{0,n}=\sum_{i=0}^{n-1} V_{i,i+1}=\sum_{i=0}^{n-1}V_{0,1}\circ \theta^i.$$
The ergodic theorem then implies that $n^{-1}V_{0,n}$ converges a.s. as $n\to\infty$. The limit must be constant
since it is a shift-invariant function of the i.i.d. sequence $(\xi_n)_{n\in\mathbb{Z}}$ defined by
$$\xi_n :=\Big((B_{n+t}-B_n)_{0\leq t\leq 1}, 
\mathfrak{R}_{[0,1]}(\nn_\infty\circ \theta^n)\Big),$$
where $\mathfrak{R}_{[0,1]}(\gamma)$ stands for the restriction of $\gamma$ to $[0,1]\times\S$. Our claim 
\eqref{voisin-tech} follows. Finally, for arbitrary $a<b$, the convergence in the lemma follows
from the special case $a=0,b=1$ using scaling and translation invariance properties of the model. \endproof

Our goal is to prove that a result similar to Lemma \ref{voisin-infty} holds for the free pointed Brownian disk.
We need a couple of preliminary lemmas.

\begin{lemma}
\label{local-bdry}
Let $\eta>0$ and $\delta>0$. Then a.s. there exists a (random) real $\ve_0>0$ such that the following holds
for every $0<\ve<\ve_0$: 
for every $u\in\mathfrak{H}_\infty^{[-\eta,\eta]}$, the property $D^\infty(u,\R)<\ve$ implies that
there exists $v\in[-\eta-\delta,\eta+\delta]$ such that $D^\infty(u,v)<\ve$, and morever
\begin{equation}
\label{local-bdry-t1}
D^\infty(u,v)= \build{\inf_{u_0=u,u_1,\ldots,u_{p-1},u_p=v}}
_{u_1,\ldots,u_{p-1}\in \mathfrak{H}_\infty^{[-\eta-2\delta,\eta+2\delta]}}^{} 
\sum_{i=1}^p D^{\infty,\circ}(u_{i-1},u_i).
\end{equation}
\end{lemma}

\proof We start by observing that
we have a.s.
\begin{equation}
\label{distech1}
\inf \Big\{D^\infty\Big(u,(-\infty,-\eta-\delta]\cup[\eta+\delta,\infty)\Big):
 u\in \mathfrak{H}_\infty^{[-\eta,\eta]}\Big\}>0
\end{equation}
and 
\begin{equation}
\label{distech2}
\inf\Big\{D^\infty\Big(u,[-\eta-\delta,\eta+\delta]\Big) 
: u\in \mathfrak{H}_\infty\backslash\mathfrak{H}_\infty^{[-\eta+2\delta,\eta+2\delta]}\Big\}>0.
\end{equation}

Let us prove \eqref{distech1}. We argue by contradiction. If the infimum in
\eqref{distech1} is zero, this means that we can find a sequence $(u_n)_{n\geq 1}$
in $\mathfrak{H}_\infty^{[-\eta,\eta]}$ such that $D^\infty(u_n,(-\infty,-\eta-\delta]\cup[\eta+\delta,\infty))$
tends to $0$ as $n\to\infty$. By compactness, we may assume that
$u_n\la u_\infty \in \mathfrak{H}_\infty^{[-\eta,\eta]}$ as $n\to\infty$
(in the sense of the topology of $\mathfrak{H}_\infty$). Then, using the continuity of the 
mapping $(u,v)\mapsto D^\infty(u,v)$, we have necessarily
$D^\infty(u_\infty,(-\infty,-\eta-\delta]\cup[\eta+\delta,\infty))=0$ and (using Lemma \ref{dist-infty})
this is only possible if there exists $v\in(-\infty,-\eta-\delta]\cup[\eta+\delta,\infty)$ such that
$D^\infty(u_\infty,v)=0$. This is a contradiction since we know that the equivalence class 
of any point of $\R$ in the quotient $\mathfrak{H}_\infty/\{D^\infty=0\}$ is a singleton.

The proof of \eqref{distech2} is similar. If the infimum in \eqref{distech2} is zero, we can find
a sequence $(v_n)_{n\geq1}$ in $\mathfrak{H}_\infty\backslash\mathfrak{H}_\infty^{[-\eta+2\delta,\eta+2\delta]}$
such that $D^\infty(v_n,[-\eta-\delta,\eta+\delta])$ tends to $0$,
and, thanks to Lemma \ref{dist-infty}, we can extract a subsequence converging to $v_\infty$.
The fact that $D^\infty(v_\infty,[-\eta-\delta,\eta+\delta])=0$ gives a contradiction. 

 We let $\ve_1$ and $\ve_2$ be the infima appearing in  formulas \eqref{distech1} and \eqref{distech2}
 respectively, and take $\ve_0=\ve_1\wedge \ve_2$.
The first assertion of the lemma follows from the definition of $\ve_1$. To get the second one, let
$u\in\mathfrak{H}_\infty^{[-\eta,\eta]}$ and $v\in[-\eta-\delta,\eta+\delta]$ such that $D^\infty(u,v)<\ve$.
Then, for every $\ve'$ such that $D^\infty(u,v)<\ve'<\ve$, we can find 
an integer $p\geq 1$ and $u_1,\ldots,u_{p-1}\in \mathfrak{H}_\infty$ such that
$$\sum_{i=1}^p D^{\infty,\circ}(u_{i-1},u_i)<\ve'$$
where $u_0=u$ and $u_p=v$. We claim that we must have 
$u_1,\ldots,u_{p-1}\in \mathfrak{H}_\infty^{[-\eta-2\delta,\eta+2\delta]}$. Indeed, if there
exists $j\in\{1,\ldots,p-1\}$ such that $u_j\in \mathfrak{H}_\infty\backslash  \mathfrak{H}_\infty^{[-\eta-2\delta,\eta+2\delta]}$,
then the bound
$$\sum_{i=j+1}^p D^{\infty,\circ}(u_{i-1},u_i)<\ve'$$
implies $D^\infty(u_j,v)<\ve'$, which contradicts the definition of $\ve_2$. \endproof

Let us now turn to the free pointed Brownian disk $\D^\bullet$. We recall the construction 
of $(\D^\bullet,D)$ in Section \ref{sec:BD-BHP} as a quotient of the space $\mathfrak{H}$
defined from a Poisson measure $\nn$ on $[0,1]\times \S$. Without loss of generality, we may
and will assume that $\nn$ is the restriction of $\nn_\infty$ to $[0,1]\times \S$. Then, 
for any $0\leq a<b<1$, the subset of $\mathfrak{H}$ defined by
$$\mathfrak{H}^{[a,b]}:=[a,b]\cup \Big(\bigcup_{j\in J, a\leq t_j\leq b} \t_{(\omega_j)}\Big)$$
is identified with the subset $\mathfrak{H}^{[a,b]}_\infty$ of $\mathfrak{H}_\infty$.
For $u,v\in \mathfrak{H}^{[a,b]}$, it will be useful to introduce the quantity 
$D^{\circ,[a,b]}(u,v)$: This quantity is defined by the same formula \eqref{pseudo-BD1} as $D^\circ(u,v)$,
except that, if one of the two intervals $[|u,v|]$ and $[|v,u|]$ of $\mathfrak{H}$ is not contained in
$\mathfrak{H}^{[a,b]}$ (this holds for at most one of the two intervals), we replace the infimum of labels on this interval by $-\infty$, or, equivalently, we 
disregard the infimum over this interval. Obviously, $D^{\circ,[a,b]}(u,v)\geq D^\circ(a,b)$.

The following lemma is then an analog of Lemma \ref{local-bdry}.

\begin{lemma}
\label{local-bdry2}
Let $\eta\in(0,1/8)$ and $\delta\in(0,1/8)$. Then a.s. there exists $\ve'_0>0$ such that the following holds
for every $0<\ve<\ve'_0$: 
for every $u\in\mathfrak{H}^{[\frac{1}{2}-\eta,\frac{1}{2}+\eta]}$, the property $D(u,[0,1])<\ve$ implies that
there exists $v\in[\frac{1}{2}-\eta-\delta,\frac{1}{2}+\eta+\delta]$ such that $D(u,v)<\ve$, and morever
\begin{equation}
\label{local-bdry-t2}
D(u,v)= \build{\inf_{u_0=u,u_1,\ldots,u_{p-1},u_p=v}}
_{u_1,\ldots,u_{p-1}\in \mathfrak{H}^{[\frac{1}{2}-\eta-2\delta,\frac{1}{2}+\eta+2\delta]}}^{} 
\sum_{i=1}^p D^{\circ,[\frac{1}{2}-\eta-2\delta,\frac{1}{2}+\eta+2\delta]}(u_{i-1},u_i).
\end{equation}
\end{lemma}

\proof The beginning of the proof  is exactly similar to that of Lemma \ref{local-bdry}, using the
obvious analogs of \eqref{distech1} and \eqref{distech2}, which hold thanks to the fact that no point of
$[0,1)$ is identified to another point of $\mathfrak{H}$ in the quotient $\mathfrak{H}/\{D=0\}$ (we leave the details to the reader). This leads
to the variant of formula \eqref{local-bdry-t2} 
where the quantities $D^{\circ,[\frac{1}{2}-\eta-2\delta,\frac{1}{2}+\eta+2\delta]}(u_{i-1},u_i)$ are replaced by 
$D^\circ(u_{i-1},u_i)$. So to get the statement of Lemma \ref{local-bdry2}, it suffices to
prove that the following claim holds for every $0<a<b<1$:
a.s. for $\ve$ small enough, if $u',v'\in \mathfrak{H}^{[a,b]}$ are such that
$D(u',[0,1])<\ve$, $D(v',[0,1])<\ve$, and
$D^\circ(u',v')<\ve$, then we have automatically $D^\circ(u',v')=D^{\circ,[a,b]}(u',v')$.

In order to prove our claim, we argue by contradiction. If the claim fails, then we can find 
a sequence $(\ve_n)_{n\geq 1}$ decreasing to $0$, and, for every $n\geq 1$, two points
$u^{(n)}$ and $v^{(n)}$ in $\mathfrak{H}^{[a,b]}$, such that:
\begin{itemize}
\item[\rm(a)] $D(u^{(n)},[0,1])<\ve_n$ and $D(v^{(n)},[0,1])<\ve_n$;
\item[\rm(b)] the interval $[|u^{(n)},v^{(n)}|]$ is not contained in $\mathfrak{H}^{[a,b]}$;
\item[\rm(c)] $\Lambda_{u^{(n)}}+\Lambda_{v^{(n)}}
-2\inf_{w\in[|u^{(n)},v^{(n)}|]}\Lambda_w <\ve_n$.
\end{itemize}
Recall the cyclic exploration $(\ee_s)_{s\in[0,\Sigma]}$ in Section \ref{sec:BD-BHP}. Let $s_n$ be as large as possible such that $\ee_{s_n}=u^{(n)}$ and similarly let $t_n$ be as small as possible
such that $\ee_{t_n}=v^{(n)}$. Because of property (b) we must have $t_n<s_n$ and 
$[|u^{(n)},v^{(n)}|]=\{\ee_r:r\in[s_n,\Sigma]\cup[0,t_n]\}$. Up to extracting a subsequence, we can assume
that $s_n\la s_\infty$ and $t_n\la t_\infty$ as $n\to \infty$. Set $u^{(\infty)}=\ee_{s_\infty}$
and $v^{(\infty)}=\ee_{t_\infty}$. By property (a) and the fact that
the equivalence class 
of any point of $[0,1)$ for the equivalence relation $\{D=0\}$ is a singleton, $u^{(\infty)}$ and $v^{(\infty)}$ must belong to $[0,1]$.
On the other hand, property (c) gives
$$\Lambda_{u^{(\infty)}}+\Lambda_{v^{(\infty)}}
-2\inf_{r\in[s_\infty,\Sigma]\cup[0,t_\infty]}\Lambda_{\ee_r} =0,$$
which implies in particular that $D(u^{(\infty)},v^{(\infty)})=0$. This means that $u^{(\infty)}=v^{(\infty)}$.
Then two cases may occur. Either $t_\infty<s_\infty$, which implies that 
$u^{(\infty)}$ is the root of one of the trees $\t_{(\omega_i)}$, but then, using the fact that $\{\ee_r:r\in[s_\infty,t_\infty]\}$
contains $[0,1]$, the last display would imply that the minimal value of $\beta$ over $[0,1]$ is attained
at the root of one of the trees $\t_{(\omega_i)}$, which does not hold a.s. Or $t_\infty=s_\infty$, but
then the last display shows that the minimal label on $\mathfrak{H}$ is attained at a point 
of $[0,1]$, which means that $v_*\in\partial\D$, contradicting \eqref{asymp-pro}. This contradiction completes the proof. \endproof

\begin{proposition}
\label{nhd-bdry}
For every $0\leq a<b\leq 1$, we have
$$\ve^{-2}\,\mathrm{Vol}\Big(\{u\in \mathfrak{H}^{[a,b]}: D(u,[0,1])<\ve\}\Big) \build{\la}_{\ve\to 0}^{}b-a,$$
in probability.
\end{proposition}

\proof We first show that the convergence in the proposition holds with 
$b-a$ replaced by $\kappa(b-a)$, and at the end of the proof we explain 
why $\kappa=1$. Thanks to the symmetries of the model (and also using the fact that $\mathfrak{H}^{[a,c]}=\mathfrak{H}^{[a,b]}\cup \mathfrak{H}^{[b,c]}$ if $0\leq a< b< c\leq1$), it is enough to consider the case $a=\frac{1}{2}-\eta$, 
$b=\frac{1}{2}+\eta$, where $\eta\in(0,1/8)$. We also fix $\delta\in(0,1/8)$. The idea is to combine
the convergence of Lemma \ref{voisin-infty} with an absolute continuity argument. Let us
introduce some notation. We set
$$V_\ve:=\mathrm{Vol}\Big(\{u\in \mathfrak{H}^{[\frac{1}{2}-\eta,\frac{1}{2}+\eta]}: D(u,[0,1])<\ve\}\Big)\ ,\ 
V^\infty_\ve:=\mathrm{Vol}\Big(\{u\in \mathfrak{H}_\infty^{[\frac{1}{2}-\eta,\frac{1}{2}+\eta]}:D^\infty(u,\R)< \ve\}\Big).$$
We also let $\bar V_\ve$ be the volume of the subset of $\mathfrak{H}$
consisting of all $u\in\mathfrak{H}^{[\frac{1}{2}-\eta,\frac{1}{2}+\eta]}$ such that there exist $v\in[\frac{1}{2}-\eta-\delta,\frac{1}{2}+\eta+\delta]$ and $u_1,\ldots,u_{p-1}\in \mathfrak{H}^{[\frac{1}{2}-\eta-2\delta,\frac{1}{2}+\eta+2\delta]}$ with the property
$$\sum_{i=1}^p D^{\circ,[\frac{1}{2}-\eta-2\delta,\frac{1}{2}+\eta+2\delta]}(u_{i-1},u_i) <\ve,$$
where $u_0=u$ and $u_p=v$. Similarly, we let $\bar V^\infty_\ve$ be the volume of the subset of $\mathfrak{H}_\infty$
consisting of all $u\in\mathfrak{H}_\infty^{[\frac{1}{2}-\eta,\frac{1}{2}+\eta]}$ such that there exist $v\in[\frac{1}{2}-\eta-\delta,\frac{1}{2}+\eta+\delta]$ and $u_1,\ldots,u_{p-1}\in \mathfrak{H}_\infty^{[\frac{1}{2}-\eta-2\delta,\frac{1}{2}+\eta+2\delta]}$ with the property
$$\sum_{i=1}^p D^{\infty,\circ}(u_{i-1},u_i) <\ve,$$
with the same convention for $u_0$ and $u_p$. 

By Lemma \ref{voisin-infty}, we know that $\ve^{-2}V^\infty_\ve$ converges in probability 
to $2\kappa\eta$. Lemma \ref{local-bdry} shows that $\P(\bar V^\infty_\ve=V^\infty_\ve)$ tends to 
$1$ as $\ve\to 0$, so that we also get that $\ve^{-2}\bar V^\infty_\ve$ converges in probability 
to $2\kappa\eta$. Now the point is that $\bar V^\infty_\ve$ is a measurable function
of $\nn$ (recall that we have assumed that $\nn$ is the restriction of $\nn_\infty$ to $[0,1]\times\S$)
and of the random path $(B_t)_{\frac{1}{2}-\eta-2\delta\leq t\leq \frac{1}{2}+\eta+2\delta}$,
whereas $\bar V_\ve$ is the {\it same} measurable function of $\nn$ and of the random path
$(\beta_t)_{\frac{1}{2}-\eta-2\delta\leq t\leq \frac{1}{2}+\eta+2\delta}$ --- here it is crucial that we use $D^{\circ,[\frac{1}{2}-\eta-2\delta,\frac{1}{2}+\eta+2\delta]}$
instead of $D^\circ$ in the 
definition of $\bar V_\ve$. The distribution of $(\beta_t)_{\frac{1}{2}-\eta-2\delta\leq t\leq \frac{1}{2}+\eta+2\delta}$
is absolutely continuous with respect to that of $(B_t)_{\frac{1}{2}-\eta-2\delta\leq t\leq \frac{1}{2}+\eta+2\delta}$,
with a Radon-Nikodym derivative bounded above by a constant $K$. It follows that the distribution 
of $\bar V_\ve$ is also absolutely continuous with respect to that of $\bar V^\infty_\ve$,
with a Radon-Nikodym derivative bounded by the same constant $K$ independent of $\ve$. Hence $\ve^{-2}\bar V_\ve$ also converges in probability 
to $2\kappa\eta$.
Then, Lemma \ref{local-bdry2} shows that $\P(\bar V_\ve=V_\ve)$ tends to 
$1$ as $\ve\to 0$, and we get that $\ve^{-2}V_\ve$ converges in probability 
to $2\kappa\eta$.

We have thus obtained the statement of the proposition, except for the value $\kappa=1$.
However, taking $a=0$ and $b=1$, we have $\ve^{-2}\mathrm{Vol}(\{x\in\D^\bullet:D(x,\partial \D^\bullet)<\ve\})
\la \kappa$ as $\ve\to\infty$, and, comparing with \cite[Proposition 2]{Disks}, we get that $\kappa=1$.
\endproof

\noindent{\it Proof of Theorem \ref{approx-unif-bdry}.}
 Thanks to Proposition \ref{nhd-bdry}, we may choose a sequence $(\ve_n)_{n\geq1}$ decreasing to $0$
such that, a.s. for every integer $N\geq 3$ and every $k\in\{0,1,\ldots,2^N-1\}$, we have
\begin{equation}
\label{approx-unif1}
\ve_n^{-2} \int_{\mathfrak{H}^{[k2^{-N},(k+1)2^{-N}]}} \mathrm{Vol}(\dd u)\,\mathbf{1}_{\{D(u, [0,1])<\ve_n\}}
\build{\la}_{n\to\infty}^{} 2^{-N}.
\end{equation}
By Lemma \ref{local-bdry2}, we also know that a.s. for every $N\geq 3$ and every $k\in\{0,1,\ldots,2^N-1\}$,
if $n$ is large enough (depending on the choice of $N$),
the conditions $u\in \mathfrak{H}^{[k2^{-N},(k+1)2^{-N}]}$ and $D(u,[0,1])<\ve_n$
imply that there exists $v\in[(k-1)2^{-N},(k+2)2^{-N}]$ such that $D(u,v)<\ve_n$
(in the case $k=0$, the notation $[-2^{-N},2\cdot2^{-N}]$ of course refers to
$[0,2\cdot2^{-N}]\cup [1-2^{-N},1]$, and similarly if $k=2^N-1$). From now on, we fix 
an element of the underlying probability space such that the preceding property and \eqref{approx-unif1} both hold.

Let $\delta>0$ and consider a bounded continuous function 
$\varphi:\D^\bullet \la \R_+$. We fix $N$ such that $|\varphi(\Pi(u))-\varphi(\Pi(v))|\leq \delta$ whenever $u,v$ both belong
to $[k2^{-N},(k+1)2^{-N}]$ for some $k\in\{0,1,\ldots,2^N-1\}$. Then, for $n$ large enough
(such that $D(x,y)<\ve_n$ implies $|\varphi(x)-\varphi(y)|<\delta$ and also such that
the property stated after \eqref{approx-unif1} holds), we have, for every $k\in\{0,1,\ldots,2^N-1\}$,
\begin{align*}
&\int_{\mathfrak{H}^{[k2^{-N},(k+1)2^{-N}]}} \mathrm{Vol}(\dd u)\,\mathbf{1}_{\{D(u, [0,1])<\ve_n\}}\,\varphi(\Pi(u))\\
&\quad\leq \Big(\max_{v\in[(k-1)2^{-N},(k+2)2^{-N}]} \varphi(\Pi(v)) +\delta\Big) \int_{\mathfrak{H}^{[k2^{-N},(k+1)2^{-N}]}} \mathrm{Vol}(\dd u)\,\mathbf{1}_{\{D(u, [0,1])<\ve_n\}}.
\end{align*}
Using  \eqref{approx-unif1}, it follows that
\begin{align*}
\limsup_{n\to\infty} \ve_n^{-2}\int_{\mathfrak{H}^{[k2^{-N},(k+1)2^{-N}]}} \mathrm{Vol}(\dd u)\,\mathbf{1}_{\{D(u, [0,1])<\ve_n\}}\,\varphi(\Pi(u))
&\leq 2^{-N}\,\Big(\max_{v\in[k2^{-N},(k+1)2^{-N}]} \varphi(\Pi(v)) +2\delta\Big)\\
&\leq \int_{k2^{-N}}^{(k+1)2^{-N}} \varphi(\Pi(v))\,\dd v + 3\delta\,2^{-N}.
\end{align*}
By summing over $k$, we get
$$\limsup_{n\to\infty} \ve_n^{-2}\int_{\mathfrak{H}} \mathrm{Vol}(\dd u)\,\mathbf{1}_{\{D(u, [0,1])<\ve_n\}}\,\varphi(\Pi(u))
\leq \int_0^1 \varphi(\Pi(v))\,\dd v +3\delta,$$
and similar arguments give the corresponding result for the liminf behavior. Since $\delta$ was arbitrary, this shows
that $\langle \mu_{\ve_n},\varphi\rangle \la \langle \mu,\varphi\rangle$, with the notation of Theorem \ref{approx-unif-bdry}. 
So we have proved that $\mu_{\ve_n}$ converges weakly to $\mu$, a.s. This is the 
desired result, except that we have restricted ourselves to a particular sequence
$(\ve_n)_{n\geq 1}$ decreasing to $0$. However, we may use \cite[Proposition 2]{Disks}, which 
already gives the a.s. weak convergence of $\mu_\ve$ as $\ve\to 0$ to a limiting
probability measure $\nu$ on $\partial \D^\bullet$. Then necessarily 
$\mu=\nu$ and this completes the proof.  \hfill$\square$

\section{Conditioning the distinguished point to belong to the boundary}
\label{sec:dis-bdry}

Our goal in this section is to give a description of the  free Brownian disk  pointed at
a point chosen uniformly on the boundary.
To this end, we will condition the free pointed Brownian disk $\D^\bullet$ on the event $\{D(v_*,\partial\D^\bullet)\leq\ve\}$
and pass to the limit $\ve\to 0$.

We first need to introduce the (non-pointed) free Brownian disk $\D$. 
To this end, write $\D^\circ$ for the (non-pointed) space obtained by forgetting the 
distinguished point of $\D^\bullet$. Then the distribution of the free Brownian disk $\D$ 
is given by the identity
$$\E[F(\D)]= \E\Big[ \frac{1}{\mathrm{vol}(\D^\circ)}\,F(\D^\circ)\Big]$$
for any nonnegative measurable function
$F$ on the space $\M$ of Section \ref{sec:GHP}
(see \cite[Section 1.5]{BM}). Conversely, we can recover the distribution 
of $\D^\bullet$ from that of $\D$ via the formula
\begin{equation}
\label{pointed-no}
\E[F( \D^\bullet)]=\E\Big[\int_\D \mathrm{Vol}(\dd x)\,F([\D,x])\Big],
\end{equation}
where we use the notation
$[\D,x]$ for the space $\D$ pointed at $x$, and $F$ is now defined on $\M^\bullet$. 
Let us give a brief justification of \eqref{pointed-no}.
For every $r>0$, write $\D_{(r)}$, resp. $\D^\bullet_{(r)}$, for the Brownian disk, resp. the pointed Brownian disk,
of perimeter $1$ and volume $r$. Note that $\D^\bullet_{(r)}$ is constructed by the very same method as
in Section \ref{sec:BD-BHP} under the probability measure $\P(\cdot\,|\,\Sigma=r)$, and that
$\D_{(r)}$ is derived from $\D^\bullet_{(r)}$ by forgetting the distinguished point. Hence, the
law of $\D^\bullet$ is obtained by integrating the law of $\D^\bullet_{(r)}$ with respect to the
density $(2\pi r^3)^{-1/2}\,\exp(-1/(2r))$ of $\Sigma$, and similarly the law of $\D$ is obtained by integrating the law of $\D_{(r)}$ with respect to the
density $(2\pi r^5)^{-1/2}\,\exp(-1/(2r))$ (see \cite[Section 1.5]{BM}).
From these considerations, we see that \eqref{pointed-no} is a consequence of the identity
\begin{equation}
\label{pointed-no2}
\E[F( \D^\bullet_{(r)})]=\E\Big[\frac{1}{r}\int_{\D_{(r)}} \mathrm{Vol}(\dd x)\,F([\D_{(r)},x])\Big].
\end{equation}
This identity follows from Lemma 18 in \cite{BM}, which is itself derived from the (trivial)
discrete analog of \eqref{pointed-no2} for quadrangulations. To be precise, \cite{BM} considers
Brownian disks as random metric spaces, without including the volume measures, but the 
argument of \cite{BM} immediately extends to our setting.

We next observe that the uniform measure $\mu$ on the boundary also makes sense
for $\D$ (it may be defined by the almost sure approximation in Theorem \ref{approx-unif-bdry}).
We then define the free Brownian disk with perimeter $1$ pointed at
a uniform boundary point as the pointed compact measure metric space $\bar \D^\bullet$, whose
distribution is given by 
$$\E[F(\bar \D^\bullet)]=\E\Big[\int_{\partial\D} \mu(\dd x)\,F([\D,x])\Big].$$

\begin{proposition}
\label{conv-pointed}
The conditional distribution of the pointed metric measure space $\D^\bullet$ given that $D(v_*,\partial\D^\bullet)\leq\ve$
converges as $\ve\to 0$ to the distribution of the free Brownian disk with perimeter $1$ pointed at
a uniform boundary point.
\end{proposition}

\proof Let $F$ be a bounded
continuous function on $\M^\bullet$, and assume that
$0\leq F\leq 1$. Then,
\begin{equation}
\label{cond-expec}
\E[F(\D^\bullet)\,|\, D(v_*,\partial\D^\bullet)\leq\ve]=
\frac{\E[F(\D^\bullet)\,\mathbf{1}_{\{D(v_*,\partial\D^\bullet)\leq\ve\}}]}{\P(D(v_*,\partial\D^\bullet)\leq\ve)},
\end{equation}
and we know from \eqref{asymp-pro} that 
$$\lim_{\ve \to 0} \ve^{-2}\,\P(D(v_*,\partial\D^\bullet)\leq\ve) = 1.$$
On the other hand, \eqref{pointed-no} gives
\begin{equation}
\label{formu-volu}
\E\Big[F(\D^\bullet)\,\mathbf{1}_{\{D(v_*,\partial\D^\bullet)\leq\ve\}}\Big]
=  \E\Big[\int_{\D} \mathrm{Vol}(\dd x)\,F([\D,x])\,\mathbf{1}_{\{D(x,\partial\D)\leq\ve\}}\Big].
\end{equation}
From Theorem \ref{approx-unif-bdry} (and the continuity of the mapping $x\mapsto [\D,x]$), we have  
$$\ve^{-2}\int_{\D} \mathrm{Vol}(\dd x)\,F([\D,x])\,\mathbf{1}_{\{D(x,\partial\D)\leq\ve\}}
\build{\la}_{\ve\to 0}^{a.s.} \int \mu(\dd y)\,F([\D,y]).$$
The desired convergence of $\E[F(\D^\bullet)| D(v_*,\partial\D^\bullet)\leq\ve]$ to $\E[\int \!\mu(\dd y)F([\D,y])]$ will follow 
from \eqref{cond-expec} and \eqref{formu-volu} if we can prove that the convergence in the last display also holds
for expected values. Arguing along a sequence of values of $\ve$
tending to $0$, Fatou's lemma gives
\begin{equation}
\label{Fatou12}
\liminf_{\ve\to 0} \ve^{-2}\E\Big[\int_{\D} \mathrm{Vol}(\dd x)\,F([\D,x])\,\mathbf{1}_{\{D(x,\partial\D^\bullet)\leq\ve\}}\Big]
\geq \E\Big[\int \mu(\dd y)\,F([\D,y])\Big].
\end{equation}
By the case $F=1$ of \eqref{formu-volu}, we have
$$\ve^{-2}\E\Big[\int_{\D} \mathrm{Vol}(\dd x)\,\mathbf{1}_{\{D(x,\partial\D)\leq\ve\}}\Big]
=\ve^{-2}\P(D(v_*,\partial\D^\bullet)\leq\ve) \build{\la}_{\ve\to 0}^{} 1.$$
Replacing $F$ by $1-F$ in \eqref{Fatou12}, we get the corresponding upper bound for the limsup, and we
conclude that we have
$$\lim_{\ve\to0} \ve^{-2}\E\Big[\int_{\D} \mathrm{Vol}(\dd x)\,F([\D,x])\,\mathbf{1}_{\{D(x,\partial\D^\bullet)\leq\ve\}}\Big]
= \E\Big[\int \mu(\dd y)\,F([\D,y])\Big].$$
This completes the proof. \endproof

We will now combine Proposition \ref{conv-bridge} and Proposition \ref{conv-pointed} to get our construction of the
free Brownian disk with perimeter $1$ pointed at
a uniform boundary point. We start from a pair $(\mathbf{b},\nn')$, where 
$\mathbf{b}=(\mathbf{b}_t)_{0\leq t\leq 1}$ is a five-dimensional Bessel bridge from $0$ 
to $0$ over the time interval $[0,1]$ and, conditionally on $\mathbf{b}$, $\nn'(\dd t\dd\omega)$ is a Poisson measure
on $[0,1]\times \S$ with intensity 
$$2\,\mathbf{1}_{\{W_*(\omega)\geq -\sqrt{3}\,\mathbf{b}_t\}}\,\dd t\,\N_0(\dd \omega).$$
 We write
$$\nn'=\sum_{j\in J} \delta_{(t'_j,\omega'_j)},$$
and $\Sigma':=\sum_{j\in J}\sigma(\omega'_j)$. 
From $\nn'$, we can define a compact measure metric space $\mathfrak{H}'$, and the 
associated cyclic exploration $(\ee'_s)_{s\in[0,\Sigma']}$, in exactly 
the same way as $\mathfrak{H}$ and $(\ee_s)_{s\in[0,\Sigma]}$ were defined from $\nn$ at the beginning of Section \ref{sec:BD-BHP}.
Intervals in $\mathfrak{H}'$ are defined as previously from the exploration $(\ee'_s)_{s\in[0,\Sigma']}$,
and we now specify labels $(\Lambda'_u)_{u\in\mathfrak{H}'}$ by setting
$\Lambda'_t:=\sqrt{3}\,\mathbf{b}_t$ for $t\in[0,1]$, and for $u\in\t_{(\omega'_j)}$, $j\in J$, 
$$\Lambda'_u:=\sqrt{3}\,\mathbf{b}_{t'_j} + \ell_u(\omega'_j).$$
A fundamental difference is now that $\Lambda'_u\geq 0$ for every $u\in\mathfrak{H}'$
(because by construction $W_*(\omega'_j)\geq -\sqrt{3}\,\mathbf{b}_{t'_j}$ for every $j\in J$).
Furthermore $0$ is the unique element of $\mathfrak{H}'$ with zero label.

We use the analogs of \eqref{pseudo-BD1} and \eqref{pseudo-BD2}, with $\Lambda_u$
replaced by $\Lambda'_u$, to define $D'^\circ(u,v)$ and  $D'(u,v)$ for $u,v\in \mathfrak{H'}$.
Then $D'(u,v)$ is a pseudo-metric on $\mathfrak{H}'$. Furthermore,
it is immediate that $D'(0,u)=\Lambda'_u$ for every $u\in \mathfrak{H}'$, and that the
bound $ |\Lambda'_u-\Lambda'_v|\leq D'(u,v)$ holds for every $u,v$. 

\begin{theorem}
\label{cons-disk}
The quotient space $\D':=\mathfrak{H}'/\{D'=0\}$ equipped with the metric induced by $D'$, with the volume 
measure which is the pushforward of the volume measure on $\mathfrak{H}'$, and with the distinguished
point which is the equivalence class of $0$, is a free Brownian disk with perimeter $1$ pointed at a uniform boundary point.\end{theorem}

We write $\Pi'$ for the canonical projection from $\mathfrak{H}'$ onto $\D'$. As previously, we can 
define the label $\Lambda'_x$ of $x\in\D'$ by setting $\Lambda'_x:=\Lambda'_u$,
for any $u\in \mathfrak{H}'$ such that $\Pi'(u)=x$. In a way similar to the formula $D(v_*,u)=\Lambda_u-\Lambda_*$
for the free pointed Brownian disk, 
labels $\Lambda'_x$ exactly correspond to distances from the distinguished point $0$ lying
on the boundary. 

\proof Thanks to Proposition \ref{conv-pointed}, it is enough to verify that the (pointed measure metric) space
$\D^\bullet$ conditioned on the event $\{D(v_*,\partial\D^\bullet)\leq \ve\}$ converges in distribution 
to $\D'$ as $\ve\to 0$,
 in the sense of the pointed Gromov-Hausdorff-Prokhorov topology. For the sake of
simplicity, we will content ourselves with proving the pointed Gromov-Hausdorff convergence
(a few additional technicalities, using Lemma 4 in \cite{Disks}, yield the stronger Gromov-Hausdorff-Prokhorov
convergence). 

Let 
$\nn$ be as previously a Poisson point measure on $[0,1]\times \S$ with
intensity $2\,\dd t\,\N_0(\dd \omega)$. Let us fix $\eta>0$ and $\alpha>0$. We first choose $\delta\in (0,1/2)$
small enough so that
\begin{equation}
\label{bridge-small}
\P\Big(\sup_{t\in[0,\delta]\cup[1-\delta,1]} \sqrt{3}\,\mathbf{b}_t <\alpha\Big)>1-\eta,
\end{equation}
and 
\begin{equation}
\label{poisson-small}
\P\Big(W^*(\omega)\leq \alpha\hbox{ for every atom }(t,\omega)\hbox{ of }\nn\hbox{ such that }t\in[0,\delta]\cup[1-\delta,1]\Big)
\geq 1-\eta.
\end{equation}
For every $\ve>0$, let $\be^\ve$ be distributed as in Proposition \ref{conv-bridge}. Thanks to 
\eqref{bridge-small} and to the
convergence in distribution in the latter proposition, we can find $\ve_0>0$ such that, for 
every $\ve\in(0,\ve_0]$, we have also
\begin{equation}
\label{bridge-small2}
\P\Big(\sup_{t\in[0,\delta]\cup[1-\delta,1]} \sqrt{3}\,\be^\ve_t <\alpha\Big)>1-\eta.
\end{equation}
Furthermore, we can also fix $\ve_1\in(0,\ve_0)$ such that, for every $\ve\in(0,\ve_1]$,
the total variation distance
between the distribution of $(\be^\ve_t+\ve)_{\delta\leq t\leq 1-\delta}$ and the 
distribution of $(\mathbf{b}_t)_{\delta\leq t\leq 1-\delta}$ is smaller than $\eta$.

Let us fix $\ve\in(0,\ve_1\wedge (\alpha/\sqrt{3})]$. On a suitable probability space, we can 
construct both $\mathbf{b}$ and $\be^\ve$ so that
\begin{equation}
\label{key-couple}
\P\Big(\ve+\be^\ve_t=\mathbf{b}_t\hbox{ for every }t\in[\delta,1-\delta]\Big) > 1-\eta.
\end{equation}
We may also assume that the Poisson point measure $\nn$ is defined on the same
probability space, and is independent of the pair $(\mathbf{b},\be^\ve)$. We then
define two other random point measures  $\nn'$ and $\nn_\ve$ by
$$\nn'(\dd t\dd \omega):= \mathbf{1}_{\{W_*(\omega)\geq -\sqrt{3}\,\mathbf{b}_t\}}\,\nn(\dd t\dd \omega)\quad
\hbox{and}\quad \nn_\ve(\dd t\dd \omega):= \mathbf{1}_{\{W_*(\omega)\geq -\sqrt{3}(\ve+\be^\ve_t)\}}\,\nn(\dd t\dd \omega).$$
Note that the pair $(\mathbf{b},\nn')$ has the same distribution as described before the
statement of the theorem, and so we may assume that the pointed metric space $\D'$
is constructed from this pair as explained above. 

Similarly, we claim that the pair $(\be^\ve,\nn_\ve)$
has the conditional distribution of the pair $(\be,\nn)$ introduced in Section \ref{sec:BD-BHP}
to construct the free pointed Brownian disk,
given the event $\{\Lambda_*\geq -\ve\sqrt{3}\}$. Let us briefly explain this. Noting
that $\{\Lambda_*\geq -\ve\sqrt{3}\}=\{W_*(\omega)\geq -\sqrt{3}(\be_{t}+\ve)\hbox{ for every atom }(t,\omega)
\hbox{ of }\nn\}$, we get, using \eqref{hittingpro},
$$\P(\Lambda_*\geq -\ve\sqrt{3}\,|\,\be)= \exp\Big(-\int_0^1 \frac{\dd t}{(\ve +\be_t)^2}\Big).$$
Then, conditionally on $\be$, the law of $\nn$ given the event $\{\Lambda_*\geq \ve\sqrt{3}\}$
is the law $\PP^{(\ve,\be)}$ of a Poisson point measure with intensity $2\,\mathbf{1}_{\{W_*(\omega)\geq-\sqrt{3}(\be_t+\ve)\}}\,
\dd t\,\N_0(\dd\omega)$. Hence, for any nonnegative measurable function $G$ on $M_p([0,1]\times\S)$,
$$\frac{\E[G(\nn)\,\mathbf{1}_{\{\Lambda_*\geq-\ve\sqrt{3}\}} \mid \be\,]}
{\P(\Lambda_*\geq -\ve\sqrt{3}\,|\,\be)} = \int \PP^{(\ve,\be)}(\dd\gamma)\,G(\gamma).$$
If $F$ is a nonnegative measurable function on $C([0,1],\R_+)$, it follows that
$$\E[F(\be)G(\nn)\,\mathbf{1}_{\{\Lambda_*\geq -\ve\sqrt{3}\}}]=\E\Big[F(\be)\exp\Big(-\int_0^1 \frac{\dd t}{(\ve +\be_t)^2}\Big)
\int \PP^{(\ve,\be)}(\dd\gamma)\,G(\gamma)\Big]= C_\ve\,\E[F(\be^\ve)\,G(\nn_\ve)],$$
which gives our claim. 

Recall the construction of the free pointed Brownian disk from the pair $(\be,\nn)$ at the beginning of Section \ref{sec:BD-BHP}. If in this construction we replace the pair $(\be,\nn)$ with $(\be^\ve,\nn_\ve)$, we can define a space $\mathfrak{H}_\ve$
analogous to $\mathfrak{H}$ and assign labels $(\Lambda^\ve_u)_{u\in\mathfrak{H}_\ve}$
to the points of $\mathfrak{H}_\ve$ (with $\Lambda^\ve_u=\sqrt{3}\,\be^\ve_t+\ell_u(\omega)$ if
$u\in \t_{(\omega)}$, for any atom $(t,\omega)$ of $\nn_\ve$). We then consider the quotient space $\D^{(\ve)}=\mathfrak{H}_\ve/\{D_\ve=0\}$,
where $D^\circ_\ve$ and then the pseudo-metric $D_\ve$ are defined by the analogs of formulas 
\eqref{pseudo-BD1} and \eqref{pseudo-BD2} in terms of the labels $\Lambda^\ve_u$. We write $v^\ve_*$ for the 
(unique) point with minimal label in $\mathfrak{H}_\ve$, and $\Pi_\ve$ for the canonical projection from $\mathfrak{H}_\ve$
onto $\D^{(\ve)}$. The preceding claim shows that $\D^{(\ve)}$ (viewed as
a compact metric space pointed at $v^\ve_*$) has the conditional
distribution of $\D^\bullet$ given the event $\{D(v_*,\partial\D^\bullet)\leq\sqrt{3}\, \ve\}$.
So, to complete the proof, we now need to a get an upper bound on the 
Gromov-Hausdorff distance between $\D^{(\ve)}$ and $\D'$. 

Let $\mathfrak{R}_\delta(\nn')$ and $\mathfrak{R}_\delta(\nn_\ve)$ denote the respective restrictions of $\nn'$
and $\nn_\ve$ to $[\delta,1-\delta]\times \S$. Similarly, let 
$\wh{\mathfrak{R}}_\delta(\nn')$ and $\wh{\mathfrak{R}}_\delta(\nn_\ve)$ denote the respective restrictions of $\nn'$
and $\nn_\ve$ to $([0,\delta)\cup(1-\delta,1])\times \S$. We consider the event $E_\ve$ where the 
following properties hold:
\begin{itemize}
\item[(i)] $\displaystyle{\sup_{t\in[0,\delta]\cup[1-\delta,1]}}\sqrt{3}(\mathbf{b}_t \vee \be^\ve_t)<\alpha$;
\item[(ii)] $\ve+\be^\ve_t=\mathbf{b}_t$ for every $t\in[\delta,1-\delta]$;
\item[(iii)] for every atom $(t,\omega)$ of $\wh{\mathfrak{R}}_\delta(\nn')$ or of $\wh{\mathfrak{R}}_\delta(\nn_\ve)$,
we have $W^*(\omega)\leq \alpha$.
\end{itemize}
Thanks to \eqref{bridge-small}, \eqref{poisson-small}, \eqref{bridge-small2} and 
\eqref{key-couple}, we have $\P(E_\ve)\geq 1-4\eta$. From now on, we argue 
on the event $E_\ve$. It is convenient to write
$$\nn_\ve=\sum_{i\in I_\ve} \delta_{(t^\ve_i,\omega^\ve_i)}$$
and to introduce the subset $\mathfrak{H}^{[\delta,1-\delta]}_\ve$ of $\mathfrak{H}_\ve$ defined by
$$\mathfrak{H}^{[\delta,1-\delta]}_\ve:=[\delta,1-\delta] \cup 
\Bigg(\bigcup_{\{i\in I_\ve: \delta\leq t^\ve_i\leq 1-\delta\}} \t_{(\omega^\ve_i)}\Bigg)$$
where we recall that the root of $\t_{(\omega^\ve_i)}$ is identified with $t^\ve_i$. 
In a similar way, we define the subset $\mathfrak{H}'^{[\delta,1-\delta]}$ of $\mathfrak{H}'$. From property
(ii) and the way $\nn_\ve$ and $\nn'$ have been constructed, we have 
 $\mathfrak{R}_\delta(\nn')=\mathfrak{R}_\delta(\nn_\ve)$, and so $\mathfrak{H}^{[\delta,1-\delta]}_\ve$
 is identified with $\mathfrak{H}'^{[\delta,1-\delta]}$.

 We define a correspondence $\r_\ve$ between the spaces
$\D^{(\ve)}$ and $\D'$ by saying that $(x,x')\in \D^{(\ve)}\times \D'$
belongs to $\r_\ve$ if and only if (at least) one of the following properties hold:
\begin{itemize} 
\item[(a)] $x=\Pi_\ve(u)$ and $x'=\Pi'(u)$, for some $u\in \mathfrak{H}^{[\delta,1-\delta]}_\ve=\mathfrak{H}'^{[\delta,1-\delta]}$;
\item[(b)] $x=\Pi_\ve(u)$ and $x'=\Pi'(0)$, for some $u\in \mathfrak{H}_\ve\backslash \mathfrak{H}^{[\delta,1-\delta]}_\ve$;
\item[(c)] $x=\Pi_\ve(v^\ve_*)$ and $x'=\Pi'(u)$, for some $u\in \mathfrak{H}'\backslash \mathfrak{H}'^{[\delta,1-\delta]}$;
\end{itemize}
We note that $\r_\ve$ contains the pair $(\Pi_\ve(v^\ve_*),\Pi'(0))$ consisting of the 
respective distinguished points of $\D^{(\ve)}$ and $\D'$.

We now need to bound the distortion of $\r_\ve$. To this end, we first observe that
$0\leq\Lambda'_u\leq 2\alpha$ if $u\in \mathfrak{H}'\backslash \mathfrak{H}'^{[\delta,1-\delta]}$ (by properties
(i) and (iii)), and similarly $-\alpha\leq \Lambda^\ve_u\leq 2\alpha$ if $u\in \mathfrak{H}_\ve\backslash \mathfrak{H}^{[\delta,1-\delta]}_\ve$. Here we use the fact that $\ve\le\alpha/\sqrt{3}$ to obtain that $\sqrt{3}\,\be^\ve_t+W_*(\omega)\geq -\alpha$ for every
atom $(t,\omega)$ of $\nn_\ve$, and we note that the lower bound $\Lambda^\ve_u\geq -\alpha$
thus holds for every $u\in\mathfrak{H}_\ve$, and in particular for $u=v^\ve_*$. Furthermore, we have 
\begin{equation}
\label{ident-label}
\ve\sqrt{3}+\Lambda^\ve_u=\Lambda'_u
\quad\hbox{ for every }u\in \mathfrak{H}^{[\delta,1-\delta]}_\ve=\mathfrak{H}'^{[\delta,1-\delta]},
\end{equation}
by property (ii) and our construction of the Poisson measures $\nn^\ve$ and $\nn'$.

For $u,v\in\mathfrak{H}_\ve$, resp. for $u,v\in \mathfrak{H}'$, let us introduce the notation $[|u,v|]^\ve$, resp. $[|u,v|]'$, for the interval 
from $u$ to $v$ in $\mathfrak{H}_\ve$, resp. in $\mathfrak{H}'$. 
For $u,v\in\mathfrak{H}_\ve$, we use the bounds 
\begin{equation}
\label{dist-bound1}
D_\ve(u,v)\leq D^\circ_\ve(u,v)=\Lambda^\ve_u+\Lambda^\ve_v -2\max\Big(\min_{w\in[|u,v|]^\ve}\Lambda^\ve_w,
\min_{w\in[|v,u|]^\ve}\Lambda^\ve_w\Big)\leq \Lambda^\ve_u+\Lambda^\ve_v+2\alpha
\end{equation}
and $D_\ve(u,v)\geq |\Lambda^\ve_u-\Lambda^\ve_v|$. Similarly, we have for $u,v\in \mathfrak{H}'$,
\begin{equation}
\label{dist-bound2}
D'(u,v)\leq D'^\circ(u,v)=\Lambda'_u+\Lambda'_v -2\max\Big(\min_{w\in[|u,v|]'}\Lambda'_w,
\min_{w\in[|v,u|]'}\Lambda'_w\Big)\leq \Lambda'_u+\Lambda'_v
\end{equation}
and $D'(u,v)\geq |\Lambda'_u-\Lambda'_v|$. We recall that $D'(0,u)=\Lambda'_u$ for every $u\in \mathfrak{H}'$,
and similarly we note that $D_\ve(v^\ve_*,u)=\Lambda^\ve_u-\Lambda^\ve_{v^\ve_*}$ for every $u\in \mathfrak{H}_\ve$.

Let $(x,x')$ and $(y,y')$ be two pairs in $\r_\ve$. In order to bound $|D_\ve(x,y)-D'(x',y')|$, we need
to distinguish several cases. Suppose first that $(x,x')$ and $(y,y')$ both satisfy property 
(b) above, so that $x'=y'=\Pi'(0)$ and there exist $u,v\in  \mathfrak{H}_\ve\backslash \mathfrak{H}^{[\delta,1-\delta]}_\ve$
such that $x=\Pi_\ve(u)$ and $y=\Pi_\ve(v)$. Then of course $D'(x',y')=0$, whereas the preceding 
bounds on labels give $D_\ve(x,y)\leq \Lambda^\ve_u+\Lambda^\ve_v+2\alpha \leq 6\alpha$, and thus 
 $|D_\ve(x,y)-D'(x',y')|\leq 6\alpha$. The same arguments show that $|D_\ve(x,y)-D'(x',y')|\leq 4\alpha$
 if $(x,x')$ and $(y,y')$ both satisfy (c). 
 
 Then suppose that $(x,x')$ satisfies (b) and $(y,y')$
 satisfies (c), and pick $u\in \mathfrak{H}_\ve\backslash \mathfrak{H}^{[\delta,1-\delta]}_\ve$ and 
 $v\in \mathfrak{H}'\backslash \mathfrak{H}'^{[\delta,1-\delta]}$ such that $x=\Pi_\ve(u)$
 and $y'=\Pi'(v)$. We have $D_\ve(x,y)=D_\ve(u,v^\ve_*)=\Lambda^\ve_u-\Lambda^\ve_{v^\ve_*}\leq 3\alpha$, whereas 
 $D'(x',y')=D'(0,v)=\Lambda'_v\leq 2\alpha$. So, in that case, we get $|D_\ve(x,y)-D'(x',y')|\leq 3\alpha$.
 
 Consider next the case where $(x,x')$ satisfies (a) and $(y,y')$ satisfies (b). Pick $u\in  \mathfrak{H}^{[\delta,1-\delta]}_\ve=\mathfrak{H}'^{[\delta,1-\delta]}$ such that $x=\Pi_\ve(u)$ and $x'=\Pi'(u)$, and $v\in \mathfrak{H}_\ve\backslash \mathfrak{H}^{[\delta,1-\delta]}_\ve$ such that $y=\Pi_\ve(v)$. Note that $\Lambda'_u=\Lambda^\ve_u + \ve\sqrt{3}$, and 
 $\ve\sqrt{3}\leq \alpha$. Since $y'=\Pi'(0)$, we have $D'(x',y')=D'(u,0)=\Lambda'_u$. On the other
 hand, $D_\ve(x,y)=D_\ve(u,v)\geq |\Lambda^\ve_u-\Lambda^\ve_v|\geq \Lambda^\ve_u-2\alpha$, and 
 $D_\ve(x,y)\leq \Lambda^\ve_u+\Lambda^\ve_v +2\alpha \leq \Lambda^\ve_u+4\alpha$. Using \eqref{ident-label}, it 
 follows that $|D_\ve(x,y)-D'(x',y')|\leq 4\alpha$. Similar arguments show that $|D_\ve(x,y)-D'(x',y')|\leq 3\alpha$
 if $(x,x')$ satisfies (a) and $(y,y')$ satisfies (c).
 
It remains to consider the more delicate case when $(x,x')$ and $(y,y')$ both satisfy (a). 
In this case, we use the following lemma.

\begin{lemma} 
\label{techGH}
On the event $E_\ve$, we have, for every $u,v\in
\mathfrak{H}^{[\delta,1-\delta]}_\ve=\mathfrak{H}'^{[\delta,1-\delta]}$,
$$|D_\ve(u,v)-D'(u,v)| \leq 7\alpha.$$
\end{lemma}

\proof
Let us fix $u,v\in
\mathfrak{H}^{[\delta,1-\delta]}_\ve=\mathfrak{H}'^{[\delta,1-\delta]}$. We observe that at least one of the following two properties holds:
\begin{itemize}
\item $[|u,v|]_\ve=[|u,v|]'$ and this interval is contained in $\mathfrak{H}^{[\delta,1-\delta]}_\ve=\mathfrak{H}'^{[\delta,1-\delta]}$;
\item $[|v,u|]_\ve=[|v,u|]'$ and this interval is contained in $\mathfrak{H}^{[\delta,1-\delta]}_\ve=\mathfrak{H}'^{[\delta,1-\delta]}$.
\end{itemize}
Without loss of generality, we may and will assume that the first property holds (otherwise, interchange 
$u$ and $v$). Note that we have then
$$\min_{w\in[|u,v|]'}\Lambda'_w=\ve\sqrt{3} + \min_{w\in[|u,v|]^\ve}\Lambda^\ve_w$$
by \eqref{ident-label}. Furthermore, the interval $[|v,u|]^\ve$ is not contained 
in $\mathfrak{H}^{[\delta,1-\delta]}_\ve$ if and only if $[|v,u|]'$ is not contained in 
$\mathfrak{H}'^{[\delta,1-\delta]}$, and then both these intervals contain $0$. 
In that case $\min_{w\in[|v,u|]'}\Lambda'_w=0$, so that the maximum appearing
in the formula for $D'^\circ(u,v)$ in \eqref{dist-bound2} must be equal to $\min_{w\in[|u,v|]'}\Lambda'_w$.

If $[|v,u|]^\ve$ is contained 
in $\mathfrak{H}^{[\delta,1-\delta]}_\ve$, the same holds for $[|v,u|]'$, and the formulas 
for $D_\ve^\circ(u,v)$ and $D'^\circ(u,v)$ in \eqref{dist-bound1} and \eqref{dist-bound2}
show that $D_\ve^\circ(u,v)=D'^\circ(u,v)$. If $[|v,u|]^\ve$ is not contained 
in $\mathfrak{H}^{[\delta,1-\delta]}_\ve$, then either the maximum 
in the formula for $D_\ve^\circ(u,v)$ in \eqref{dist-bound1} is $\min_{w\in[|u,v|]^\ve}\Lambda^\ve_w$,
and this implies again that $D_\ve^\circ(u,v)=D'^\circ(u,v)$, or this maximum 
is $\min_{w\in[|v,u|]^\ve}\Lambda^\ve_w$, which belongs to $[-\alpha,0]$, and we 
have
\begin{equation}
\label{dist-bound3}
D^\circ_\ve(u,v) \geq \Lambda^\ve_u + \Lambda^\ve_v.
\end{equation}

Next, from the definition of $D_\ve(u,v)$ as an infimum, we can find an integer $p\geq 1$
and $u_0,u_1,\ldots,u_p\in \mathfrak{H}_\ve$ such that $u_0=u$, $u_p=v$, and
$$\sum_{i=1}^p D_\ve^\circ(u_{i-1},u_i) \leq D_\ve(u,v)+\alpha.$$
We then distinguish several cases:
\begin{itemize}
\item All $u_i$, $i\in\{1,\ldots,p-1\}$, belong to $\mathfrak{H}^{[\delta,1-\delta]}_\ve=\mathfrak{H}'^{[\delta,1-\delta]}$, and, for every 
$j\in\{1,\ldots,p\}$, the maximum in the formula for $D_\ve^\circ(u_{j-1},u_j)$
is attained for one of the two intervals $[|u_{j-1},u_j|]^\ve$ and $[|u_{j},u_{j-1}|]^\ve$ 
that is contained in $\mathfrak{H}^{[\delta,1-\delta]}_\ve$. In that case, the preceding considerations
show that we have $D_\ve^\circ(u_{j-1},u_j)=D'^\circ(u_{j-1},u_j)$, for every 
$j\in\{1,\ldots,p\}$, and then
$$D'(u,v)\leq \sum_{i=1}^p D'^\circ(u_{i-1},u_i)=\sum_{i=1}^p D_\ve^\circ(u_{i-1},u_i)\leq D_\ve(u,v)+\alpha.$$
\item All $u_i$, $i\in\{1,\ldots,p-1\}$, belong to $\mathfrak{H}^{[\delta,1-\delta]}_\ve=\mathfrak{H}'^{[\delta,1-\delta]}$, but, for some
$j\in\{1,\ldots,p\}$, the maximum in the formula for $D_\ve^\circ(u_{j-1},u_j)$
is attained for one of the two intervals $[|u_{j-1},u_j|]^\ve$ and $[|u_{j},u_{j-1}|]^\ve$ 
that is not contained in $\mathfrak{H}^{[\delta,1-\delta]}_\ve$. In that case, using \eqref{dist-bound3} together with the
lower bound $D^\circ_\ve(u_{i-1},u_i)\geq |\Lambda^\ve_{u_i}-\Lambda^\ve_{u_{i-1}}|$, and then
\eqref{ident-label}, we have
\begin{align*}
D_\ve(u,v)\geq -\alpha + \sum_{i=1}^p D_\ve^\circ(u_{i-1},u_i)
&\geq -\alpha + \sum_{i=1}^{j-1} |\Lambda^\ve_{u_i}-\Lambda^\ve_{u_{i-1}}| + \Lambda^\ve_{u_{j-1}}+ \Lambda^\ve_{u_j}
+  \sum_{i=j+1}^{p} |\Lambda^\ve_{u_i}-\Lambda^\ve_{u_{i-1}}|\\
&\geq -\alpha +\Lambda^\ve_u + \Lambda^\ve_v\\
&=-\alpha + \Lambda'_u + \Lambda'_v -2\ve\sqrt{3}\\
&\geq D'(u,v)-3\alpha.
\end{align*}
\item For some $j\in\{1,\ldots,p-1\}$, $u_j\notin \mathfrak{H}^{[\delta,1-\delta]}_\ve$,
which implies $\Lambda^\ve_{u_j}\leq 2\alpha$. In that case,
\begin{align*}
D_\ve(u,v)\geq -\alpha + \sum_{i=1}^p D_\ve^\circ(u_{i-1},u_i)
\geq -\alpha + |\Lambda^\ve_u-\Lambda^\ve_{u_j}| + |\Lambda^\ve_v-\Lambda^\ve_{u_j}|
&\geq -\alpha + \Lambda^\ve_u + \Lambda^\ve_v -4\alpha\\
&= \Lambda'_u + \Lambda'_v - 2\ve\sqrt{3}- 5\alpha\\
&\geq D'(u,v)-7\alpha.
\end{align*}
\end{itemize}
So in all cases we have obtained $D_\ve(u,v)\geq D'(u,v)-7\alpha$. 
A symmetric argument (in fact simpler since we do not have to consider
the second case) shows similarly that $D'(u,v)\geq D_\ve(u,v)-7\alpha$. 
This completes the proof of Lemma \ref{techGH}. \endproof

Let us complete the proof of Theorem \ref{cons-disk}. 
If the pairs $(x,x')$ and $(y,y')$ in $\r_\ve$ both satisfy property (a) above, we
can find $u,v\in
\mathfrak{H}^{[\delta,1-\delta]}_\ve=\mathfrak{H}'^{[\delta,1-\delta]}$ such that
$x=\Pi_\ve(u),x'=\Pi'(u)$ and $y=\Pi_\ve(v),y'=\Pi'(v)$. Then Lemma \ref{techGH} gives
$$|D_\ve(x,y)-D'(x',y')|=|D_\ve(u,v)-D'(u,v)|\leq 7\alpha.$$
Recalling the bounds obtained before the statement of Lemma \ref{techGH},
we conclude that the distortion of $\r_\ve$ is bounded above
by $7\alpha$. 

From the relation between the Gromov-Hausdorff distance and the infimum
of distortions of correspondences (\cite[Proposition 3.6]{LGM}, see \cite[Theorem 7.3.25]{BBI} for a proof in the
non-pointed case, which is easily adapted), we get that, on the event $E_\ve$, 
the (pointed) Gromov-Hausdorff distance between $\D^{(\ve)}$
and $\D'$ is bounded above by $7\alpha/2$. Since $E_\ve$
has probability at least $1-4\eta$, and both $\alpha$ and $\eta$
were arbitrary, this shows that $\D^{(\ve)}$ converges in probability
to $\D'$,  in the sense of the pointed Gromov-Hausdorff topology.
This completes the proof. \endproof

\section{The Caraceni-Curien construction of the Brownian half-plane}
\label{sec:CC}

In this section, we show that the construction of the Brownian half-plane from \cite{BMR,GM0},
which is presented in Section \ref{sec:BD-BHP}, is equivalent to
the construction proposed by Caraceni and Curien \cite[Section 5.3]{CC}.
This is a direct application of our construction of the free Brownian disk 
pointed at a uniform boundary point (Theorem \ref{cons-disk}).

We start by recalling the construction of \cite{CC}. We consider a random process 
$(X_t)_{t\in\R}$ such that $(X_t)_{t\geq 0}$ and $(X_{-t})_{t\geq 0}$
are two independent five-dimensional Bessel processes started from $0$.
We then consider a random point measure $\nn^\star(\dd t\dd \omega)$
on $\R\times \S$, such that, conditionally on $X$, $\nn^\star$ is Poisson
with intensity
$$2\,\mathbf{1}_{\{W_*(\omega)\geq-\sqrt{3}\,X_t\}}\,\dd t\,\N_0(\dd \omega).$$
We write 
$$\nn^\star=\sum_{i\in I^*} \delta_{(t_i^\star,\omega_i^\star)}$$
and, in a way similar to the previous sections, we let $\mathfrak{H}^\star$ be obtained from the disjoint union
$$\R\cup \Bigg(\bigcup_{i\in I^\star} \t_{(\omega^\star_i)}\Bigg),$$
by identifying the root of $\t_{(\omega^\star_i)}$ with $t_i^\star$, for 
every $i\in I^\star$. As in the preceding section, for any reals $a<b$, we introduce the subset of 
$\mathfrak{H}^\star$ defined by
$$\mathfrak{H}^{\star,[a,b]}:=[a,b]\cup \Bigg(\bigcup_{i\in I^\star, t_i\in[a,b]} \t_{(\omega^\star_i)}\Bigg).$$

We assign nonnegative labels to the points of $\mathfrak{H}^\star$
by setting $\Lambda^\star_u:=\sqrt{3}\,X_u$ if $u\in \R$, and $\Lambda^\star_u
:= \sqrt{3}\,X_{t_i^\star} + \ell_u(\omega_i^\star)$ if $u\in  \t_{(\omega^\star_i)}$, $i\in I^\star$.
The exploration process $(\ee^*_t)_{t\in\R}$ of $\mathfrak{H}^\star$ is then defined 
as in Section \ref{sec:BD-BHP}, and allows us to consider intervals $[|u,v|]$, for 
$u,v\in\mathfrak{H}^\star$. The functions $D^{\star,\circ}(u,v)$ and $D^\star(u,v)$,
for $u,v\in\mathfrak{H}^\star$, are 
defined by the exact analogs of formulas \eqref{pseudo-infty1} and \eqref{pseudo-infty2},
just replacing $\Lambda^\infty$ with $\Lambda^\star$. Then, it is straightforward
to verify that $D^\star(0,u)=D^{\star,\circ}(0,u)=\Lambda^\star_u$, for every $u\in \mathfrak{H}^\star$.

We set $\H^\star:=\mathfrak{H}^\star/\{D^\star=0\}$, and equip $\H^\star$ with the 
distance induced by $D^\star$, and with the distinguished point that is the equivalence class of $0$.
We view $(\H^\star,D^\star)$ as a random pointed boundedly compact length space (we could also
consider the volume measure on $\H^\star$, but we refrain from doing so for the sake of
simplicity). The compactness of closed balls in $\H^\star$ follows from the property
\begin{equation}
\label{transience-lab}
\lim_{a\to\infty} \inf\{\Lambda^\star_u:u\in \mathfrak{H}^\star\backslash \mathfrak{H}^{\star,[-a,a]}\}=+\infty,
\end{equation}
whose easy proof is left to the reader (see \cite[Lemma 3.3]{CLG} for a very similar argument).

We note that $(\H^\star,D^\star)$ is scale invariant, meaning that, for every 
$\lambda>0$, the space obtained when
multiplying the metric $D^\star$ by the factor $\lambda$ has the same distribution
as the original space. This scale invariance property is a straightforward
consequence of our construction.

\begin{theorem}
\label{ident-CC}
The pointed boundedly compact length spaces $(\H,D^\infty)$ and $(\H^\star,D^\star)$ have the same
distribution.
\end{theorem}

\proof Recall the space $(\D',D')$ in Theorem \ref{cons-disk}, which is a free Brownian
disk with perimeter $1$ pointed at a uniform boundary point. For every 
$r>0$, we write $B_r(\D')$ for the closed ball 
of radius $r$ centered at the distinguished point of $(\D',D')$, and we use 
the similar notation $B_r(\H)$ or $B_r(\H^\star)$. We view these balls as 
(random) pointed compact metric spaces. For $\lambda>0$, we also use the notation 
$\lambda\cdot B_r(\D')$ for the ``same'' space with the metric multiplied by
the factor $\lambda$.

 It follows from
Corollary 3.9 in \cite{BMR} that $(\H,D^\infty)$ is the tangent cone in distribution of
$(\D',D')$ at its distinguished point, meaning that, for every $r>0$, $\lambda\cdot B_{r/\lambda}(\D')$ converges 
in distribution to $B_r(\H)$ when $\lambda\to\infty$, in the sense of the pointed Gromov-Hausdorff topology. So, to get the theorem, it will be enough to
prove that $\lambda\cdot B_{r/\lambda}(\D')$ converges in distribution to $B_r(\H^\star)$ as $\lambda\to\infty$.
Recalling the scale invariance property of $(\H^\star,D^\star)$, this follows immediately 
from the next lemma, which is analogous to \cite[Corollary 3.9]{BMR}.

\begin{lemma}
\label{CC-tech}
For every $\delta>0$, we can find $\ve_0>0$ such that, for every $0<\ve\leq \ve_0$, 
we can couple the spaces $(\D',D')$ and $(\H^\star,D^\star)$ in such a way that the balls
$B_\ve(\D')$ and $B_\ve(\H^\star)$ are equal with probability at least $1-\delta$.
\end{lemma}

\proof We fix $\delta>0$. Recall the notation introduced before Theorem \ref{cons-disk}, and the construction of $(\D',D')$ as a quotient space of the 
(labeled) space $\mathfrak{H}'$, which is defined from the pair $(\mathbf{b},\nn')$. For every $\eta\in(0,1/2)$, let $\mathfrak{H}'^{(\eta)}$ be the subset
of $\mathfrak{H}'$ defined by
$$\mathfrak{H}'^{(\eta)}:=[0,\eta] \cup[1-\eta,1]\cup \Bigg(\bigcup_{j\in J, t'_j\in[0,\eta] \cup[1-\eta,1]} \t_{(\omega'_j)}\Bigg),$$
with the same identifications as previously ($0$ is identified to $1$, and the root of $\t_{(\omega'_j)}$
is identified to $t'_j$).

By choosing $\eta\in(0,1/2)$ sufficiently small,
we may couple the pairs $(\mathbf{b},\nn')$ and $(X,\nn^*)$ in such a way that the following
two properties hold with probability at least $1-\delta/2$:
\begin{itemize}
\item[(i)] $\mathbf{b}_t=X_t$ and $\mathbf{b}_{1-t}=X_{-t}$ for every $t\in[0,\eta]$;
\item[(ii)] the restriction of $\nn'(\dd t\dd \omega)$ to $[0,\eta]\times \S$ coincides with the restriction
of $\nn^\star(\dd t\dd \omega)$ to $[0,\eta]\times \S$, and the restriction of $\nn'(\dd t\dd \omega)$ to $[1-\eta,1]\times \S$ coincides with the pushforward of the restriction
of $\nn^\star(\dd t\dd \omega)$ to $[-\eta,0]\times \S$ under the translation $(t,\omega)\mapsto (1+t,\omega)$.
\end{itemize}
In what follows, we assume that $(\mathbf{b},\nn')$ and $(X,\nn^*)$ have been coupled in this way. 
Then, we can choose $\ve>0$ small enough so that the two properties
\begin{itemize}
\item[(iii)] $\Lambda'_u>3\ve$ for every $u\in \mathfrak{H}'\backslash \mathfrak{H}'^{(\eta)}$,
\item[(iv)] $\Lambda^\star_u>3\ve$ for every $u\in \mathfrak{H}^\star\backslash \mathfrak{H}^{\star,[-\eta,\eta]}$,
\end{itemize}
hold except on a set of probability at most $\delta/2$ (we use \eqref{transience-lab} for (iv)).

Let $E_{\eta,\ve}$ be the event where 
properties (i)---(iv) hold, so that $\P(E_{\eta,\ve})\geq 1-\delta$. We will verify that
the property $B_\ve(\D')=B_\ve(\H^\star)$ holds on $E_{\eta,\ve}$. From now on we argue on the event $E_{\eta,\ve}$. 

We first observe that there is an obvious one-to-one correspondence between the sets 
$\mathfrak{H}'^{(\eta)}$ and $\mathfrak{H}^{\star,[-\eta,\eta]}$. In particular, the subset
$[1-\eta,1]$ of $\mathfrak{H}'^{(\eta)}$ corresponds to the subset $[-\eta,0]$ of $\mathfrak{H}^{\star,[-\eta,\eta]}$
via the translation $u\mapsto u-1$ and (thanks to (ii)) the trees $\t_{(\omega'_j)}$ for indices $j$ such that
$1-\eta<t'_j<1$ correspond to the trees $\t_{(\omega^\star_i)}$ for indices $i$ such that $-\eta<t^\star_i<0$.
If $u\in\mathfrak{H}'^{(\eta)}$, we will write $u^\star$ for the corresponding point in $\mathfrak{H}^{\star,[-\eta,\eta]}$.
Using (i) and (ii), we have $\Lambda'_u=\Lambda^\star_{u^\star}$ for every $u\in\mathfrak{H}'^{(\eta)}$. 
Moreover, if $u,v\in \mathfrak{H}'^{(\eta)}$, the set $[|u^\star,v^\star|]\cap \mathfrak{H}^{\star,[-\eta,\eta]}$
coincides with the set of all points $w^\star$ for $w\in [|u, v|]\cap \mathfrak{H}'^{(\eta)}$. 

Recall that $\Lambda'_u=D'(0,u)$ for every $u\in\mathfrak{H}'$, and $\Lambda^\star_v=D^\star(0,v)$
for every $v\in \mathfrak{H}^\star$. It follows from (iii) and (iv) that the sets $\{u\in\mathfrak{H}':D'(0,u)\leq 3\ve\}$
and $\{v\in \mathfrak{H}^\star:D^\star(0,v)\leq3\ve\}$ are contained in $\mathfrak{H}'^{(\eta)}$ and 
in $\mathfrak{H}^{\star,[-\eta,\eta]}$ respectively, and these two sets are equal modulo the preceding
identification of $\mathfrak{H}'^{(\eta)}$ with $\mathfrak{H}^{\star,[-\eta,\eta]}$. To complete the proof,
it remains to verify that, for any $u,v\in \mathfrak{H}'$ such that $D'(0,u)\leq \ve$ and $D'(0,v)\leq \ve$, we have $D'(u,v)=D^\star(u^*,v^*)$. 

To get this, first observe that, if $ \wt u, \wt v\in \mathfrak{H}'$ are such that $D'(0, \wt u)\leq 3\ve$ and $D'(0,\wt v)\leq 3\ve$
(so that in particular $\wt u,\wt v \in \mathfrak{H}'^{(\eta)}$), we have
$$\inf_{w\in [|\wt u,\wt v|]}\Lambda'_w=\inf_{w\in [|\wt u, \wt v|]\cap \mathfrak{H}'^{(\eta)}}\Lambda'_w
=\inf_{w\in [|\wt u^\star,\wt v^\star|]\cap \mathfrak{H}^{\star,[-\eta,\eta]}} \Lambda^\star_w= 
\inf_{w\in [|\wt u^\star,\wt v^\star|]} \Lambda^\star_w,$$
and thus $D'^\circ( \wt u, \wt v)=D^{\star,\circ}(\wt u^\star,\wt v^\star)$. 

Then, let $u,v\in \mathfrak{H}'$ such that $D'(0,u)\leq \ve$ and $D'(0,v)\leq \ve$. We observe that, in the definition 
of $D'(u,v)$ as an infimum
over choices of $u_0=u,u_1,\ldots,u_{p-1},u_p=v$, we may disregard the case when  $D'(0,u_j)=\Lambda'_{u_j}>3\ve$
for some $j\in\{1,\ldots,p-1\}$, because, if this happens,  
the lower bound $D'^\circ(u_{i-1},u_i)\geq |\Lambda'_{u_i}-\Lambda'_{u_{i-1}}|$ implies that
$$\sum_{i=1}^p D'^\circ(u_{i-1},u_i) \geq |\Lambda'_{u_j}-\Lambda'_u| +|\Lambda'_{u_j}-\Lambda'_v|\geq 4\ve,$$
whereas we have $D'(u,v)\leq D'(0,u)+D'(0,v)\leq 2\ve$.
A similar remark applies to the definition of $D^\star(u^\star,v^\star)$, and we conclude from the 
equality $D'^\circ( \wt u, \wt v)=D^{\star,\circ}(\wt u^\star,\wt v^\star)$ when 
$D'(0, \wt u)\leq 3\ve$ and $D'(0,\wt v)\leq 3\ve$ that the infima in the respective definitions of $D'(u,v)$ and $D^\star(u^*,v^*)$ are equal, so that 
$D'(u,v)=D^\star(u^*,v^*)$ as desired. This completes the proof of Lemma \ref{CC-tech} and Theorem \ref{ident-CC}. 
\endproof


\begin{thebibliography}{99}

\bibitem{ADH}
{R. Abraham, J.-F. Delmas, P. Hoscheit},
A note on Gromov-Hausdorff-Prokhorov distance between (locally) compact measure spaces.
{\it Electron. J. Probab.} 18, 1--21 (2013)

\bibitem{ALG}
{C. Abraham, J.-F. Le Gall}, Excursion theory for Brownian motion
indexed by the Brownian tree. {\it J. Eur. Math. Soc. (JEMS)} 20, 2951--3016 (2018)

\bibitem{AHS}
{M. Albenque, N. Holden, X. Sun}, Scaling limit of large triangulations of polygons. Preprint, available
at arXiv:1910.04946.


\bibitem{BMR}
E. Baur, G. Miermont, G. Ray,
Classification of scaling limits of uniform quadrangulations with a boundary.
{\it Ann. Probab.} 47, 3397--3477 (2019) 

\bibitem{Bet0}
{J. Bettinelli}, Scaling limits for random quadrangulations of positive genus.
Electron. J. Probab. 15, 1594--1644 (2010)

\bibitem{Bet}
{J. Bettinelli}, Scaling limit of random planar 
quadrangulations with a boundary. {\it Ann. Inst. H. Poincar\'e Probab. Stat.} 51, 432--477 (2015)

\bibitem{BM}
{J. Bettinelli, G. Miermont},
Compact Brownian surfaces I. Brownian disks. 
{\it Probab. Theory Related Fields} 167, 555-614 (2017)

\bibitem{BR}
{T. Budzinski, A. Riera},
Personal communication.


\bibitem{BBI}
{D. Burago, Y. Burago, S. Ivanov}, {\it A Course in Metric Geometry}.
Graduate Studies in Mathematics, vol. 33. Amer. Math. Soc., Boston, 2001.

\bibitem{CC}
{A. Caraceni, N. Curien},
Geometry of the Uniform Infinite Half-Planar Quadrangulation.
{\it Random Structures Algorithms} 52, 454--494 (2018)

\bibitem{Plane}
{N. Curien, J.-F. Le Gall}, 
The Brownian plane. 
{\it J. Theoret. Probab.} 27, 1240--1291 (2014)


\bibitem{CLG} N. Curien, J.-F. Le Gall, The hull process of the Brownian plane.
{\it Probab. Theory Related Fields} 166, 187--231 (2016)

\bibitem{Getoor}
R.K. Getoor, The Brownian escape process. {\it Ann. Probab.} 7, 864--867 (1979)


\bibitem{GM0}
{E. Gwynne, J. Miller},
Scaling limit of the uniform infinite half-plane quadrangulation in the Gromov-Hausdorff-Prokhorov-uniform topology.
{\it Electron. J. Probab.} 22, paper no 84, 47pp. (2017)


\bibitem{GM2}
{E. Gwynne, J. Miller},
Convergence of the free Boltzmann quadrangulation with simple boundary to the Brownian disk.
{\it Ann. Inst. Henri Poincar\'e Probab. Stat.} 55, 1--60 (2019)

\bibitem{Zurich} {J.-F. Le Gall}, {\it Spatial Branching Processes, Random Snakes and 
Partial Differential Equations}. {Lectures in Mathematics ETH Z\"urich}. Birkh\"auser, Boston, 1999.

\bibitem{Uniqueness} {J.-F. Le Gall}, Uniqueness and universality of the Brownian map.
{\it Ann. Probab.} {41}, 2880--2960 (2013)


\bibitem{LGM}
{J.-F. Le Gall, G. Miermont},
Scaling limits of random trees and planar maps. In: Probability and Statistical Physics in Two and More Dimensions, 
Clay Mathematics Proceedings, vol.15, pp.155-211, AMS-CMI, 2012


\bibitem{Disks} J.-F. Le Gall, Brownian disks and the Brownian snake. {\it Ann. Inst. H. Poincar\'e Probab. Stat.} 55, 237--313 (2019)

\bibitem{spine}
{J.-F. Le Gall, A. Riera},
Spine representations for non-compact  models of random geometry.
Preprint, 2020.

\bibitem{Miller}
{J. Miller}, Liouville quantum gravity as a metric space and a scaling limit.
Proceedings of ICM 2018, available at arXiv:1712.01571


\bibitem{MS}
{J. Miller, S. Sheffield},  Liouville quantum gravity and the Brownian map II: geodesics and continuity of the embedding.
Preprint, available at arXiv:1605.03563 

\bibitem{PY2}
{J.W. Pitman, M. Yor},
Bessel processes and infinitely divisible laws.
In: Stochastic Integrals. Lecture Notes Math. 851, pp. 285--370. Springer 1981.


\bibitem{PY}
{J.W. Pitman, M. Yor},
A decomposition of Bessel bridges. {\it Z. Wahrsch. verw. Gebiete} 59, 425--457 (1982)

\bibitem{RY} {D. Revuz, M. Yor}, {\it
Continuous Martingales and Brownian
Motion.} Springer, Berlin, 1991.

\bibitem{Vervaat}
{W. Vervaat}, A relation between Brownian bridge and Brownian excursion.
{\it Ann. Probab.} 7, 143--149 (1979)

\end{thebibliography}
\end{document}